\documentclass[11pt]{amsart}
\usepackage[all]{xy}
\usepackage{amssymb}

\newtheoremstyle{obs}% name
  {3pt}%      Space above
  {3pt}%      Space below
  {}%         Body font
  {}%         Indent amount (empty = no indent, \parindent = para indent)
  {\bfseries}% Thm head font
  {.}%        Punctuation after thm head
  {.5em}%     Space after thm head: " " = normal interword space;
  {}%         Thm head spec (can be left empty, meaning `normal')
\theoremstyle{obs}
\numberwithin{equation}{section}

\newtheorem{definition}{Definition}[section]

\newtheorem{proposition}[definition]{Proposition}

\newtheorem{remarkth}[definition]{Remark}

\newenvironment{remark}{\begin{remarkth}\upshape}{\hfill$\diamond$\end{remarkth}}
\renewcommand{\emph}[1]{{\bfseries\itshape{#1}}}
\usepackage{amssymb}

\newcommand{\lcf}{\lbrack\! \lbrack} %Corchete de algebroide de Lie izquierda
\newcommand{\rcf}{\rbrack\! \rbrack} %Corchete de algebroide de Lie derecha.

\begin{document}

\title[Generalized variational calculus ]{\textbf{\textsc{Generalized variational calculus for continuous and discrete mechanical systems}}}

\author[V. D{\'\i}az]{Viviana D{\'\i}az}
\address{Viviana D{\'\i}az: Departamento de Matem\'atica, Universidad Nacional del Sur,
8000 Bah{\'\i}a Blanca,
Argentina } \email{viviana.diaz@uns.edu.ar}

\author[D.\ Mart\'{\i}n de Diego]{David Mart\'{\i}n de Diego}
\address{David Mart\'{\i}n de Diego:
Instituto de Ciencias Matem\'aticas %\linebreak
(CSIC-UAM-UC3M-UCM),
c$\backslash$ Nicol\'as Cabrera, 13-15, Campus Cantoblanco,UAM
28049 Madrid, Spain}\email{david.martin@icmat.es}

\date{\today}

\thanks{This work has been partially supported by UNS, Argentina (project PGI 24/ZL06); FONCYT,
Argentina (project PICT 2010-2746); CONICET, Argentina (project PIP 2010--2012 11220090101018); MEC (Spain) Grants MTM2013-42870-P, MTM2009-08166-E, and IRSES-project ``Geomech-246981''.}

\maketitle

%\tableofcontents
%\begin{index}
%\printtheindex
\begin{abstract}
In this paper, we consider a generalization of variational calculus which allows us to consider in the same framework different cases of mechanical systems, for ins\-tan\-ce, Lagrangian mechanics, Hamiltonian mechanics, systems subjected to constraints, optimal control theory and so on. This generalized variational calculus is based on two main notions: the tangent lift of curves and the notion of complete lift of a vector field. Both concepts are also adapted for the case of skew-symmetric algebroids, therefore, our formalism easily extends to the case of Lie algebroids and nonholonomic systems (see also \cite{GraGra08}). Hence, this framework automatically includes reduced mechanical systems subjected or not to constraints. Finally, we show that our formalism can be used  to tackle the case of discrete mechanics, including reduced systems, systems subjected to cons\-traints and discrete optimal control theory.
\end{abstract}

\section{\textbf{\textsc{ Introduction}}}

The main objective of classical mechanics is to seek for trajectories describing the motion
of mechanical systems and its properties.  It is well-known that there exists a variational procedure to obtain these
trajectories for many cases of interest.  Hamilton's variational principle singles out particular curves $q: [t_0, t_1]\rightarrow {\mathbb R}$ by
\[
\delta\int_{t_0}^{t_1} L(q(t), \dot{q}(t))\; {\rm d}t=0\; ,
\]
where the variation is over curves   joining two
fixed points.
A basic result of  calculus
of variations is that Hamilton's variational principle (see \cite{foundation}) holds for a curve $q(t)$ if and only if the curve satisfies the Euler-Lagrange equations:
\[
\frac{d}{dt}\left(\frac{\partial L}{\partial \dot{q}}\right)- \frac{\partial L}{\partial {q}}=0\; .
\]
The variational derivation of the equations of motion are extended to many systems of interest; for instance, in the dynamic of systems
 associated with Lie
groups, one can derive  the Euler-Poincar\'e equations which occur for many
systems; e.g., rigid body equations, equations of fluids and
plasma dynamics \cite{Holm1,Holm2}. For other systems, as
an spacecraft with movable internal parts, one can combine Euler-Poincar\'e and Euler-Lagrange equations, both derive from appropriate variational procedures.

In this paper, we explore the common features of all these systems obtaining  a ge\-ne\-ra\-li\-zed variational derivation of the equations of motion.
Our method is  valid for a wide class of mechanical systems including Lagrangian and Hamiltonian mechanics, variational systems with constraints,
nonholonomic systems and reduced systems. Moreover, the techniques are easily adapted for the case of discrete mechanics.  More specifically,
we define a generalized variational problem on $TQ$ only  determining a submanifold $\Sigma$ of $T^*TQ$ where $Q$ stands for the configuration
space of a mechanical system. Then, using the notions of tangent lift of curves  and vector fields  (see Section \ref{sect:GVCtangbundle} for more details), we extend Hamilton's variational principle in the following way:
a solution of a generalized variational problem determined by $\Sigma\subset T^*TQ$ is a curve
$\sigma:I\rightarrow Q$ such that
\begin{equation}\label{pou}
\int_I\left\langle\mu(t),X^T(t,\dot{\sigma}(t))\right\rangle \,{\rm d}t=0,
\end{equation}
where $\mu$ is a curve in the submanifold $\Sigma$ which projects over $\sigma,$ and $X^T$ is the tangent  lift to $TQ$ of an arbitrary  time-dependent vector field on $Q$.

%\end{enumerate}
We will show that these generalized variational problems accomplishe  a great number of systems of interest in mechanics.
Additionally, since our approach is intrinsic, we may derive the corresponding  Hamel's formalism where the velocity components are
measured relative to a set of independent vector fields on the configuration space $Q$  not ge\-ne\-ra\-lly associated
with  configuration coordinates.
Moreover, it is possible to substitute the tangent bundle by another space which admits the lifting operations ne\-ce\-ssa\-ry for our de\-fi\-ni\-tion of generalized variational calculus. One example of this type of spaces is precisely skew-symmetric algebroids which allows us to define the co\-rres\-pon\-ding equations of motion. With the general framework of skew-symmetric algebroids, we derive the equations for interesting type of mechanical systems: Euler-Poincar\'e equations, Lie-Poisson equations, Lagrange-Poincar\'e equations, equations for nonholonomic systems, higher-order lagrangian mechanics and so on.
These applications for  continuous  lagrangian systems were studied previously in \cite{GraGra08} where the authors develop a variational calculus adapted to skew-symmetric algebroids, finding the equations for lagrangian systems in this setting and also for the case of systems subjected to different type of constraints (nonholonomic or vakonomic). In this paper, we analyze the underlying geo\-me\-try of infinitesimal variational calculus allowing new and in\-te\-res\-ting applications as, for instance, discrete mechanics.
Moreover, our formalism follows the same philosophy of the classical approach to variational calculus using exterior differential systems, i.e., Griffiths formalism, in which it is given a subbundle $I$ of the cotangent bundle $T^*M$ of a manifold $M$ and a 1-form $\varphi$ on $M$. The subbundle $I$ determines the set
curves $\sigma: I\rightarrow M$ such that $\sigma^*(I)=0$ (integral curves of $I$)
and the formalism studies the extremals of the  functional $J(\sigma)=\int_{\sigma}\varphi$
(see also \cite{MR684663, MR1189496}).

In the case of discrete mechanics, we will start with
a submanifold $\Sigma_d$ of $T^*Q\times T^*Q\equiv T^*(Q\times Q)$ and, using an appropriate  discrete  tangent lift of  vector fields
(see Section \ref{DGVCenQxQ} for more details)
and discrete curves, we extend the
discrete Hamilton's variational principle (see \cite{marsden-west}). In this extension, we consider as  solutions of the discrete  generalized variational problem determined by $\Sigma_d$, the discrete curves
 $\sigma:\mathbb{Z}\rightarrow Q$ such that there exists a curve $\mu:\mathbb{Z}\rightarrow\Sigma_d\subset T^*(Q\times Q)$ which projects over the curve $\tilde{\sigma}(k)=(\sigma(k),\sigma(k+1))\in Q\times Q$ and, for all $\mathbb{Z}$-dependent section $X:\mathbb{Z}\times Q\rightarrow TQ,$
$$\sum_{k=0}^{N-1}\left\langle\mu(q_k,q_{k+1}),X^T(k,q_k,q_{k+1})\right\rangle=0$$ holds.

We will see that this description is flexible enough to cover the most important cases of discrete variational calculus, also with constraints, and even to be defined on Lie groupoids (see \cite{MMM06Grupoides,weinstein96} and references therein).

For a better understanding of our methods, we will start with the two more familiar cases of tangent bundles; namely, the continuous case and the cartesian product of two copies of the configuration space (the discrete setting). Then, we will move to the case of mechanics on skew-symmetric algebroids and Lie groupoids, showing that the techniques are quite similar to the standard cases.

\vspace{.5cm}
The paper is structured as follows. %Some elementary concepts of tangent bundle geometry are reviewed in  Section \ref{sect:fibtangente}.
\tableofcontents

\section{\textbf{\textsc{Generalized variational calculus on the tangent bundle}}}\label{sect:GVCtangbundle}

\subsection{Tangent bundle geometry}\label{sect:fibtangente}

Given a differentiable manifold $Q$ and a fixed point $x\in Q$, we can introduce the notion of \emph{curve at $x$} as a curve $\gamma: I \longrightarrow Q$  such that $I\subseteq {\mathbb R}$ contains 0 in its interior and $\gamma(0)=x$.  Then, we say that two curves $\gamma_1$ and $\gamma_2$ at $x$ are \emph{equivalent} if, for any coordinate chart $(U, \varphi)$ with $x\in U,$ we have that $\frac{d\gamma_1}{dt}(0)=\frac{d\gamma_2}{dt}(0)$. Therefore, with this definition, it is possible to introduce  an equivalence relation of curves at $x$ and define
a \emph{tangent vector} $v_x$ as an equivalence class $v_x=[\gamma]^{(1)}_x$.
The collection of all equivalence classes defines the \emph{tangent space} $T_xQ$. The \emph{tangent bundle} is precisely the disjoint union of tangent spaces
$TQ=\bigsqcup_{x\in M} T_xQ$ equipped with a natural structure of vector bundle. We denote by $\tau_{TQ}: TQ\rightarrow Q$ the
canonical projection\footnote{For a vector bundle $E$ over $Q$ we use the notation $\tau_E:E\rightarrow Q$ and $\pi_E:E^*\rightarrow Q$ for
the vector bundle projections of $E$ and $E^*.$ This is not the typical notation in the case of the tangent bundle, but we will use it for coherence with Section \ref{sect:algebroides}.} defined by $\tau_{TQ}(v_x)=x$.
Coordinates $(q^i)$ in $Q,$ for $1\leq i\leq n$ if $n=\dim Q,$ induce natural coordinates $(q^i, \dot{q}^i)$ in $TQ$ such that
$\tau_{TQ}(q^i, \dot{q}^i)=(q^i)$.  Its dual vector bundle is the \emph{cotangent bundle} $T^*Q$ with projection $\pi_{TQ}: T^*Q\rightarrow Q$ (for more details,  see \cite{foundation, ManuelMecanica}).

Similarly, it is also possible to define the \emph{second-order tangent bundle} $T^{(2)}Q$ taking  equivalence classes  of curves $\gamma_1$
and $\gamma_2$ at $x$ where $\frac{d\gamma_1}{dt}(0)=\frac{d\gamma_2}{dt}(0)$ and
$\frac{d^2\gamma_1}{dt^2}(0)=\frac{d^2\gamma_2}{dt^2}(0)$. In general, one can define higher-order tangent bundles using this procedure, see
\cite{ManuelCampos}. We alternatively denote by $[\gamma]^{(2)}_x$ or $a_x$ the corresponding equivalence class in $T^{(2)}Q$. We have induced coordinates
$(q^i, \dot{q}^i, \ddot{q}^i)$ in $T^{(2)}Q$.
In this case, we consider the canonical immersion
$j_2: T^{(2)}Q\rightarrow TTQ$ defined as
$j_2([\gamma]_x^{(2)})=[{\gamma}^{(1)}]_{v_x}^{(1)}$, where
${\gamma}^{(1)}$ is the lift of the curve $\gamma$ to
$TQ$ and ${\gamma}^{(1)}(0)=v_x$; that is,  the curve ${\gamma}^{(1)}: {\mathbb R}\rightarrow
TQ$ is given by $\gamma^{(1)}(t)=[\gamma_t]_{x}^{(1)}$
where $\gamma_t(s)=\gamma(t+s)$. In local coordinates
\[
j_2(q^{i},\dot{q}^{i}, \ddot{q}^{i})=(q^{i},\dot{q}^{i}, \dot{q}^{i}, \ddot{q}^{i})\; .
\]

Given a map $f: Q_1\rightarrow Q_2$ between two manifolds, we have the \emph{tangent map} $Tf\equiv f_*: TQ_1\rightarrow TQ_2$ defined by
\[
f_*([\gamma]^{(1)}_x)=[f\circ \gamma]^{(1)}_{f(x)}\;.
\]
 Based on this tangent lift of a map, there exists a canonical lift of a curve on $Q$  to a curve on the tangent bundle $TQ$. In fact, if we have a  curve
 $\sigma: I\rightarrow Q,$  we define the \textbf{tangent lift} of $\sigma$ as $\dot{\sigma}\equiv \frac{d\sigma}{dt}: I\rightarrow  TQ$ such that
\[
\dot\sigma(t)=T\sigma(t, 1)\in T_{\sigma(t)}Q.
\]
In coordinates, if $\sigma(t)=(q^i(t)),$ then $\dot{\sigma}(t)=(q^i(t),\dot{q}^i(t))$.

Another important geometric ingredient that we will need for our definition of ge\-ne\-ra\-li\-zed variational calculus is the  notion of complete lift of a vector field. Remember that a vector field $X$ is a smooth section of $\tau_{TQ}: TQ\rightarrow Q$, that is, $X\in \Gamma (\tau_{TQ})\equiv {\mathfrak X}(Q)$.
Expressed in terms of the coordinate frame $\{\partial /\partial q^i\}$, we have that
\[
X=X^i(q)\frac{\partial}{\partial q^i}\; .
\]

 We denote by $\{\Phi^X_t\}$ the flow of $X$.  The most natural definition of the \emph{complete lift}  $X^c$ of $X$  is given in terms of its flow.
 We say that $X^c$ is the vector field on $TQ$ with flow $\{T\Phi^X_t\}$. In other words,
\[
X^c(v_x)=\left.\frac{d}{dt}\right|_{t=0}\left(T_x\Phi^X_t(v_x)\right).
\]
%in other words, if $v_x=[\gamma]_x$ then
%\[
%X^c(v_x)=\frac{d}{dt}\huge|_{t=0}\frac{d}{ds}\huge|_{t=0}\left(\Phi^X_t(\gamma(s))\right)
%\]
In the standard coordinate frame $\{\partial/\partial q^i,\partial/\partial \dot{q}^i \}$, we have that
\begin{equation}\label{complete}
X^c=X^i(q)\frac{\partial}{\partial q^i}+\dot{q}^j\frac{\partial X^i}{\partial q^j}\frac{\partial}{\partial \dot{q}^i}.
\end{equation}

Schematically,
\begin{equation*}
\xymatrix{
TQ\ar[d]_{\tau_{TQ}}\ar[r]^{X^c}&TTQ \ar[d]^{T\tau_{TQ}\equiv (\tau_{TQ})_*}\\
Q\ar[r]^{X}&TQ
}
\end{equation*}

In our approach, we will need an alternative characterization of the complete lift. Then, recall first that a \emph{linear function} on the
vector bundle $\tau_{TQ}: TQ\rightarrow Q$ is identified to a section of the dual bundle $\pi_{TQ}:T^*Q\rightarrow Q$. More precisely, if
$\beta\in \Gamma(\pi_{TQ})$ (that is, $\beta$ is a 1-form), then
we define the linear function $\hat{\beta}: TQ\rightarrow {\mathbb R}$ by
\[
\hat{\beta}(v_x)=\langle \beta(x), v_x\rangle,
\]
for all $v_x\in TQ$.
Then, an alternative characterization of the complete lift will be the following.

\begin{proposition}\label{prop:levantcpleto}
The complete lift $X^c$ of a vector field on $Q$ is the unique vector field on $TQ$ such that
verifies the following two conditions:
\begin{enumerate}
\item $X^c$ is projectable over $X$ by means of $(\tau_{TQ})_*$, that is, $(\tau_{TQ})_*X^c=X$.
\item $X^c(\hat{\alpha})=\widehat{{\mathcal L}_X\alpha}$, for all $\alpha\in \Gamma(\pi_{TQ})$.
\end{enumerate}
\end{proposition}
Here, ${\mathcal L}_X\alpha\in \Gamma (\pi_{TQ})$ denotes the Lie derivative of $\alpha$ with respect to $X$, that is,
\[
\langle {\mathcal L}_X\alpha, Y\rangle={\mathcal L}_X \langle \alpha, Y\rangle - \langle \alpha,  [X,Y]\rangle,\  \forall\ Y\in \Gamma(\tau_{TQ}).
\]

An interesting remark is about the choice of a frame to locally write the complete lift of a vector field. In (\ref{complete}), we have used the
standard frame but, in some cases, it is interesting to use a different one. Let us assume that we
have fixed coordinates $(q^i)$ in $Q$ and an arbitrary frame $\{Y_i\}$ (a nonholonomic or moving frame, following different authors) where
\[
Y_j=\rho_j^i(q)\frac{\partial}{\partial q^i}.
\]

Then, a vector field $X\in {\mathfrak X}(Q)$ has the following local expressions
\[
X=\tilde{X}^j Y_j=\tilde{X}^j\rho_j^i\frac{\partial}{\partial q^i}.
\]
Moreover, the new frame induces a new system of coordinates $(q^i, y^i)$ on $TQ$, where $v_x=y^i Y_i(x)$ for any $v_x\in TQ$.
Using Proposition \ref{prop:levantcpleto} or by a change of coordinates, it is not hard to prove that the complete lift $X^c$ can be rewritten as
\[
X^c=\tilde{X}^j\rho^i_j\frac{\partial}{\partial q^i}+\left(\rho^i_j\frac{\partial \tilde{X}^k}{\partial q^i}- \mathcal{C}_{ij}^k
\tilde{X}^i \right) y^j\frac{\partial}{\partial y^k},
\]
where the structure function $\mathcal{C}_{ij}^k$ are defined by $[Y_i,Y_j]=\mathcal{C}_{ij}^k Y_k.$

Another notion that will be used later is the \emph{vertical lift} of a vector field on $Q$ to $TQ$. Let $X\in {\mathfrak X}(Q)$, the vertical lift of $X$
is the vector field on $TQ$ defined by:
\[
X^v(v_x)=\left.\frac{d}{dt}\right|_{t=0}(v_x+tX(x))),\; \ \forall\ v_x\in T_xQ.
\]
Locally,
\[
X^v=X^i(q)\frac{\partial}{\partial\dot{q}^i}
\]
or, in the  frame $\{Y_i\},$ we have that
\[
X^v=\tilde{X}^j(q)\frac{\partial}{\partial y^j}\, .
\]
An alternative definition of vertical lift is the following:
\begin{proposition}
   The vertical lift $X^v$ of a vector field $X$ is the unique vector field on $TQ$ verifying the following conditions:
   \begin{enumerate}
      \item $X^v(\tau^*_{TQ} f)=0,$ for all $f\in C^{\infty}(Q),$
      \item $X^v(\hat\alpha)=\tau^*_{TQ}(\langle\alpha,X\rangle),$ for all $\alpha\in\Gamma(\pi_{TQ}).$
   \end{enumerate}
\end{proposition}

\vspace{.3cm}
For our study we need to deal with time-dependent vector fields and the notion of their tangent lifts.

A \emph{time-dependent vector field} $X$ is a smooth mapping
$
X: I\times Q\rightarrow TQ,$ for $I\subseteq {\mathbb R},$ such that $X(t, x)\in T_x Q$.
We denote the set of time-dependent vector fields by ${\mathfrak X}(pr_Q)$ where
$pr_Q: I\times Q\rightarrow Q$.
%Observe that, with this notation, ${\mathfrak X}(id_Q)={\mathfrak X}(Q)$.

\begin{definition}
The \emph{tangent lift} $X^T$ of a time-dependent vector field $X$ on $Q$ is the unique time-dependent vector field on $TQ$ verifying the following two conditions:
\begin{enumerate}
\item $X^T$ is projectable over $X$ by means of $\tau_{TQ}$, that is, $(\tau_{TQ})_*X^T=X$.
\item $(X^T)_{t,v_x}(\hat{\alpha})=\widehat{{\mathcal L}_{X_t}\alpha}(v_x)+\frac{d}{dt}\langle \alpha_x, X_x(t)\rangle$, for all $\alpha\in \Gamma(\pi_{TQ})$.
\end{enumerate}
\end{definition}

Here, $X(t,x)=X_t(x)=X_x(t)$.

Schematically,
\begin{equation*}
\xymatrix{
\mathbb{R}\times TQ\ar[d]_{(id_{\mathbb{R}},\tau_{TQ})}\ar[rr]^{X^T}&&TTQ\ar[d]^{T\tau_{TQ}}\\
\mathbb{R}\times Q\ar[rr]^{X}&&TQ
}
\end{equation*}

In the coordinate frame $\{\partial/\partial q^i\}$, we have
\[
X=X^i(t,q)\frac{\partial}{\partial q^i}
\]
then
\[
X^T=X^i(t,q)\frac{\partial}{\partial q^i}+\left(\frac{\partial X^i}{\partial t}+\dot{q}^j\frac{\partial X^i}{\partial q^j}\right)\frac{\partial}{\partial \dot{q}^i}.
\]
In the frame $\{Y_i\}$ where
\[
X=\tilde{X}^i(t,q)Y_i,
\]
by using the coordinates $(q^i, y^i)$ in $TQ$ induced by the  frame $\{Y_i\},$ we have
\[
X^T=\tilde{X}^j(t,q)\rho^{i}_j(q)\frac{\partial}{\partial q^i}+\left[\frac{\partial\tilde{X}^k}
{\partial t}(t,q)+\left( \rho^{i}_j(q)\frac{\partial \tilde{X}^k}{\partial q^i}(t,q)- \mathcal{C}_{ij}^k(q)\tilde{X}^i(t, q)\right) y^j\right]\frac{\partial}{\partial y^k}.
\]

Similarly, we can introduce the \emph{vertical lift} $X^V$ of a time-dependent vector field  $X\in \mathfrak{X}(pr_Q)$ as
\[
X^V(t, v_x)=(X_t)^v(v_x),
\]
where $X_t$ is the
vector field on $Q$ defined by $X_t(x) = X(t, x)$.

In canonical coordinates $\displaystyle X^V=X^i(t,q)\frac{\partial}{\partial q^i}$ or, in the nonholonomic frame, $\displaystyle X^V=\tilde{X}^j(t,q)\frac{\partial}{\partial y^j}.$

\vspace{.2cm}
Also, we define the \emph{total derivative} of a function $f: {\mathbb R}\times Q\rightarrow {\mathbb R}$ as the function
$
\frac{df}{dt}: {\mathbb R}\times TQ\rightarrow {\mathbb R}$  defined by
\[
\frac{df}{dt}(t, v_x)=\frac{\partial f}{\partial t}(t, v_x)+v_x(f_t), \hbox{ where } (t, v_x)\in {\mathbb R}\times T_xQ\; .
\]
Locally, we have that $\displaystyle{\frac{df}{dt}=\frac{\partial f}{\partial t}+\frac{\partial f}{\partial q^i}\dot{q}^i}$.

 In the same way, if $F: {\mathbb R}\times TQ\rightarrow {\mathbb R}$, its total derivative  is the function
$
\frac{dF}{dt}: {\mathbb R}\times T^{(2)}Q\rightarrow {\mathbb R}$  defined by
\[
\frac{dF}{dt}(t, [\gamma]^{(2)}_x)=\frac{\partial F}{\partial t}(t, [\gamma]^{(1)}_x)+ j_2([\gamma]^{(2)}_x)(F_t), \hbox{ where } (t, [\gamma]^{(2)}_x)\in {\mathbb R}\times T^{(2}_xQ\; .
\]
Locally, we can write $\displaystyle{\frac{dF}{dt}=\frac{\partial F}{\partial t}+\frac{\partial F}{\partial q^i}\dot{q}^i+\frac{\partial F}{\partial \dot{q}^i}\ddot{q}^i}$.

\

The following definition will  play an important role in the sequel.

\begin{definition}\label{def:E-L-Oper}
The \emph{Euler-Lagrange operator} associated with a 1-form $\mu\in \Gamma(\pi_{TQ})=\Lambda^1(TQ)$  is the mapping ${\mathcal E}_{\mu}: T^{(2)}Q\rightarrow T^*Q$ defined by
\[
 \left\langle {\mathcal E}_{\mu} ([\gamma]^{(2)}_x), X(x)\right\rangle= \frac{d}{dt}\langle \mu,  X^v\rangle ([\gamma]^{(2)}_x)-\langle \mu, X^c\rangle ([\gamma]^{(1)}_x),
\]
 for any $X\in {\mathfrak X}(Q)$.
\end{definition}

This is well defined since the definition of the Euler-Lagrange operator only depends on the point $X(x)=v_x$.

Observe that if $X\in {\mathfrak X}(pr_Q),$ we have that
\begin{equation}\label{pok}
 \left\langle {\mathcal E}_{\mu} ([\gamma]^{(2)}_x), X(t,x)\right\rangle= \frac{d}{dt}\langle \mu,  X^V\rangle(t, [\gamma]^{(2)}_x)-\langle \mu, X^T\rangle(t, [\gamma]^{(1)}_x).
 \end{equation}

Locally, in the frame $\displaystyle\{\partial/\partial q^i\},$ if $\mu=\mu_i\mathrm{d}q^i+\tilde{\mu}_i\mathrm{d}\dot{q}^i,$ we have that
 \begin{equation}
 \displaystyle\left\langle\mathcal{E}_{\mu},\frac{\partial}{\partial q^i}\right\rangle=\displaystyle \frac{d}{dt}\tilde{\mu}_i-\mu_i
 \end{equation}
 or, in a non canonical frame $\{Y_i\}$ with coordinates $(q^i,y^i),$
 \begin{equation}\label{esubmunoncan}
 \displaystyle\left\langle{\mathcal E}_{\mu},Y_i\right\rangle =\displaystyle\frac{d}{dt}\tilde{\mu}_i-\mu_j\rho_i^j+\mathcal{C}_{ij}^k y^j\tilde{\mu}_k\; ,
 \end{equation}
where $\mu=\mu_i\mathrm{d}q^i+\tilde{\mu}_i\mathrm{d}y^i$.

 %%%%%%%%%%%%%%%%%%%%%%%%%%%%%%%%%%%%%%%%%%%%%%%

 \vspace{.3cm}
For a function $L: TQ\rightarrow \mathbb{R},$
 \[
 \left\langle {\mathcal E}_{dL} , \frac{\partial}{\partial q^i}\right\rangle= \frac{d}{dt}\left(\frac{\partial L}{\partial \dot{q}^i}\right)
 -\frac{\partial L}{\partial q^i}\;
 \]
or, in an arbitrary frame $\{Y_i\}$, for an element $X=\tilde{X}^iY_i\in {\mathfrak X}(pr_Q)$ we have that
\begin{equation}\label{EulerLagnonholframe}
 \langle {\mathcal E}_{\mathrm{d}L} , X\rangle= \tilde{X}^i\left[
 \frac{d}{dt}\left(\frac{\partial L}{\partial y^i}\right)
 -\rho^j_i\frac{\partial L}{\partial q^j}+{\mathcal C}_{ij}^k
y^j\frac{\partial L}{\partial y^k}\right]\; .
\end{equation}

\subsection{Generalized variational problem on the tangent bundle}

\begin{definition}
A \emph{generalized variational problem} on $TQ$ is determined by a submanifold $\Sigma$ of $T^*TQ$.
\end{definition}
 We initially assume the submanifold property for simplicity since in general $\Sigma$ could be any subset of $T^*TQ$.

\begin{definition}\label{solution-m}
A  \emph{solution} of the generalized variational problem determined by $\Sigma\subset T^*TQ$ is a  smooth curve
$\sigma:I\rightarrow Q$ such that
%\begin{enumerate}
%\item  $\text{Im}(\dot{\sigma})\subseteq C,$
%\item
there exists another curve $\mu:I\rightarrow\Sigma$  verifying $\pi_{TTQ}(\mu(t))=\dot{\sigma}(t)$ and, for all time-dependent vector field $X\in\mathfrak{X}(pr_Q),$
\begin{equation}\label{pou}
\int_I\left\langle\mu(t),X^T(t,\dot{\sigma}(t))\right\rangle \,\mathrm{d}t=0.
\end{equation}
%\end{enumerate}
\end{definition}

Taking the canonical projection $\pi_{TTQ}: T^*TQ\rightarrow TQ,$ we define the subset $C= \pi_{TTQ}(\Sigma)$
(\emph{kinematical constraints}) and we have schematically the following diagram

\begin{equation*}
\xymatrix{
&\Sigma\ar[d]_{(\pi_{TTQ})_{|\Sigma}}\ar@{^{(}->}[rr]^{i_{\Sigma}}&&T^*TQ\ar[d]^{\pi_{TTQ}}\\
&C\ar@{^{(}->}[rr]&&TQ\ar[d]_{\tau_{TQ}}\\
I\ar[rrr]_{\sigma}\ar[ur]^{\dot{\sigma}}\ar@/^3.5pc/[uur]^{\mu} & &&Q
}
\end{equation*}

Since $\pi_{TTQ}(\mu(t))=\dot{\sigma}(t)$ then $\text{Im}(\dot{\sigma})\subseteq C.$

\vspace{.2cm}
Applying the definition of the Euler-Lagrange operator introduced in definition \ref{def:E-L-Oper}, we deduce from equation (\ref{pou}) that
\begin{equation*}
0=\int_I
  \left(\frac{d}{dt}\langle \mu, X^V\rangle  (t, \ddot{\sigma}(t))-\langle{\mathcal E}_{\mu} (\ddot{\sigma}(t)), X(t, {\sigma}(t))\rangle)\right)\, \mathrm{d}t.
\end{equation*}

\vspace{.3cm}
If $I=[t_0,t_1],$ then $$\int_{t_0}^{t_1}\left(\frac{d}{dt}\langle \mu, X^V\rangle  (t, \ddot{\sigma}(t))-\langle{\mathcal E}_{\mu} (\ddot{\sigma}(t)), X(t, {\sigma}(t))\rangle)\right)\mathrm{d}t=0$$
is equivalent to
\begin{equation}\label{IntEsubmuX}
\int_{t_0}^{t_1}\langle{\mathcal E}_{\mu} (\ddot{\sigma}(t)), X(t,{\sigma}(t))\rangle)\mathrm{d}t=\langle \mu, X^V\rangle (t,\dot{\sigma}(t))\Big\vert_{t_0}^{t_1}.
\end{equation}

Assuming for simplicity that $X(t_0,\sigma(t_0))=X(t_1,\sigma(t_1))=0,$ we avoid the boundary conditions and, therefore, applying the
Fundamental Lemma of Calculus of Variations, we have that $\mu:I\rightarrow\Sigma\subseteq T^*TQ$ verifies the following equations

\begin{equation} \label{trescondiciones}
\left\{
\begin{array}{rcl}
  \mathcal{E}_{\mu}(\ddot{\sigma}(t))&=&0\\
   \\
    \pi_{TTQ}(\mu(t))&=&\dot\sigma(t)
\end{array}
\right..
\end{equation}

In particular, $Im(\dot{\sigma}(t))\subseteq  C.$

In canonical coordinates, if we assume that $\Sigma$ is determined by the vanishing of cons\-traints $\Phi^{\alpha}=0$ in $T^*TQ,$ a
curve $\sigma: t\rightarrow(q^i(t))$ is a solution of the generalized variational problem if there exists a 1-form $\mu$ along
$\dot{\sigma}(t);$ that is, $\mu=(\mu_ i(q,\dot{q})\mathrm{d}q^i+\tilde{\mu}_i(q,\dot{q})\mathrm{d}\dot{q}^i)\vert_{\dot{\sigma}(t)},$
such that
$$\mathcal{E}_{\mu}(q(t),\dot{q}(t),\ddot{q}(t))=0,$$ or equivalently $$\frac{\partial\tilde{\mu}_i}{\partial\dot{q}^j}\ddot{q}^j+
\frac{\partial\tilde{\mu}_i}{\partial q^j}\dot{q}^j-\mu_i=0.$$

Therefore, we write locally equations (\ref{trescondiciones}) as follows
\begin{equation}\label{rty}
\left\{
\begin{array}{rl}
   \displaystyle\frac{d}{dt}(\tilde{\mu}_i(q(t),\dot{q}(t)))-\mu_i(q(t),\dot{q}(t))&=0 \\
   \\
   \Phi^{\alpha}(q^i(t),\dot{q}^i(t),\mu_i(q(t),\dot{q}(t)),\tilde{\mu}_i(q(t),\dot{q}(t)))&=0
\end{array}
\right..
\end{equation}
It is generically difficult to obtain useful characterizations of equations (\ref{rty}), but we will see in the next subsections that for particular choices of $\Sigma$, we will derive the equations of motion of many mechanical systems of interest.

\subsection{Lagrangian mechanics}

 Given a Lagrangian function $L:TQ\rightarrow\mathbb{R}$, we know that the classical Euler-Lagrange equations for $L$ are derived using variational principles (see for instance \cite{foundation}). Of course, our generalized variational calculus is equivalent to the classical derivation using standard variational techniques.
In this particular case, we have that
$\Sigma=Im(\mathrm{d}L)=\mathrm{d}L(TQ)$ and $C=TQ.$ Observe that $\Sigma$ is a Lagrangian submanifold of $T^*TQ$ equipped with the canonical symplectic 2-form $\omega_{TQ}$.
So we look for
 a curve $\sigma:I=[t_0,t_1]\rightarrow Q$ such that
$\displaystyle\int_{t_0}^{t_1}\left\langle \mathrm{d}L(\dot{\sigma}(t)),X^T(t, \dot{\sigma}(t))\right\rangle\, \mathrm{d}t=0,$ for all $X\in {\mathfrak X}(pr_Q),$
and we also assume that $X(t_0,\sigma(t_0))=X(t_1,\sigma(t_1))=0.$

In this case, $\mu(t)=\mathrm{d}L(\dot{\sigma}(t)).$

Using Equation (\ref{pok}) we deduce that
\begin{eqnarray*}
 \displaystyle 0 & =& \int_{t_0}^{t_1}\left\langle\mathrm{d}L(\dot{\sigma}(t)),X^T(t, \dot{\sigma}(t))\right\rangle\ \mathrm{d}t\\
  &=&\int_{t_0}^{t_1}
  \left(\frac{d}{dt}\langle\mathrm{d}L,X^V\rangle(t, \ddot{\sigma}(t))-\langle{\mathcal E}_{\mathrm{d}L} (\ddot{\sigma}(t)), X(t, {\sigma}(t))\rangle)\right)\, \mathrm{d}t\\
  &=&-\int_{t_0}^{t_1}
  \langle{\mathcal E}_{\mathrm{d}L} (\ddot{\sigma}(t)), X(t, {\sigma}(t))\rangle\, dt + \left.\langle \mathrm{d}L,X^V\rangle(t, \dot{\sigma}(t))\right|_{t_0}^{t_1}.
  \end{eqnarray*}
Therefore, the equations of motion of Lagrangian mechanics are
\[
{\mathcal E}_{\mathrm{d}L}=0.
\]

Locally, in the coordinate frame, we obtain the classical Euler-Lagrange equations
\[\displaystyle\frac{d}{dt}\left(\frac{\partial L}{\partial\dot{q}^i} \right)-\frac{\partial L}{\partial q^i}=0.
\]
In the  frame $\{Y_i\},\ Y_i\in\mathfrak{X}(Q)$ for $1\leq i\leq n,$ we derive another representation of the Euler-Lagrange equations: the Hamel equations (see equation \ref{IntEsubmuX})
\begin{equation*}
\left\{
\begin{array}{rl}
\displaystyle\frac{d}{dt}\left(\frac{\partial\tilde{L}}{\partial y^i}\right)
 -\rho^j_i\frac{\partial\tilde{L}}{\partial q^j}+{\mathcal C}_{ij}^k
y^j\frac{\partial\tilde{L}}{\partial y^k}&=0 \\
\\
\dot{q}^i&=\rho^i_j y^j
\end{array}
\right.,
\end{equation*}
where $\tilde{L}(q^i,y^i)=L(q^i,\rho^i_j(q)y^j)$ and $Y_j=\rho^i_j(q)\frac{\partial}{\partial q^i}$.

\subsection{Hamiltonian mechanics}\label{mecHam}

Let $H:T^*Q\rightarrow\mathbb{R}$ be a Hamiltonian function. We will show that the typical Hamilton equations for $H$ are also expressed as a generalized
va\-ria\-tio\-nal problem. First,  we will use the canonical  antisymplectomorphism $\mathcal{R}$  between $(T^*T^*Q,\omega_{T^*Q})$ and $(T^*TQ,\omega_{TQ})$ (see
references  \cite{GraGra08,MackPX94} and references therein), that in local coordinates is given by $\mathcal{R}(q,p,\mu_q,\mu_p)=(q,\mu_p,-\mu_q,p).$

Taking  the submanifold $\mathrm{d}H(T^*Q)=Im(\mathrm{d}H)$ of  $T^*T^*Q$ and using $\mathcal{R},$ we construct the submanifold
$\Sigma_H=\mathcal{R}(\mathrm{d}H(T^*Q))$ of $T^*TQ$.
In local coordinates we can write
$$\Sigma_H=\left\{(q^i,\dot{q}^i,\mu_i,\tilde{\mu}_i) \; \mid \; \dot{q}^i =\frac{\partial H}{\partial p_i}(q, \tilde{\mu}), \ {\mu}_i=-\frac{\partial H}{\partial q^i}(q, \tilde{\mu})\right\}.$$

Given such a $\Sigma_H,$ we have the following definition.

\begin{definition}
   A curve $\sigma:I\rightarrow Q$ is a solution of the Hamiltonian problem determined by $H: T^*Q\rightarrow {\mathbb R}$ if there exists  a curve $\mu:I\rightarrow
   \Sigma_H\subset T^*TQ$ such that $\pi_{TTQ}(\mu(t))=\dot{\sigma}(t)$ and, for all $X\in\mathfrak{X}(pr_Q),$
   $$\int_I\left\langle\mu_{\dot{\sigma}(t)},X^T(t,\sigma(t))\right\rangle\mathrm{d}t=0.$$
\end{definition}

Locally, the curve $\mu:I\rightarrow\Sigma_H$
is such that $\mu:t\mapsto(q^i(t),\dot{q}^i(t),\mu_i(t),\tilde\mu_i(t))$ where $\displaystyle\dot{q}^i(t) =\frac{\partial H}{\partial p_i}(q(t), \tilde{\mu}(t))$,  $\ {\mu}_i(t)=-\frac{\partial H}{\partial q^i}(q(t), \tilde{\mu}(t)).$ Therefore, the equations of motion derived from  $\Sigma_H$
are:
\[
\displaystyle\mathcal{E}_{\mu}=\frac{d\tilde{\mu}^i}{dt}+\frac{\partial H}{\partial q^i}=0,
\]
and the equation  $Im(\dot{\sigma}(t))\subseteq C$ is now rewritten as
$\displaystyle\dot{q}^i(t) =\frac{\partial H}{\partial p_i}(q(t), \tilde{\mu}(t))$. Both equations are the typical Hamilton's equations for $H: T^*Q\rightarrow {\mathbb R}$.

\subsection{Constrained variational calculus}\label{calcvariacrestringido}

In this secton, we study the case of variational constrained calculus, also called vakonomic mechanics (see references
\cite{blochcrouch93,cardinfavretti,CdDdLM03,GraGra08,zampieri}). The equations are derived using  purely variational techniques. We will see how to define a submanifold of $T^*TQ$ to reproduce these classical equations using the generalized variational calculus.

From a geometrical point of view, these type of variationally constrained problems are determined by a pair $(C, l)$ where $C$ is a submanifold of $TQ$, with
inclusion $i_C: C\hookrightarrow TQ$, and $l: C\rightarrow \mathbb{R}$ a Lagrangian function defined \emph{only} along $C$. So we can define
 \begin{multline*}
    \Sigma _l = \bigl\{ \mu \in T ^\ast TQ \mid \pi _{TTQ} (\mu) \in C \text{
        and } \left\langle \mu, v \right\rangle = \left\langle
        \mathrm{d} l , v \right\rangle, \\
      \text{ for all } v \in T C \subset T TQ \text{ such that } \tau
      _{TTQ} (v) = \pi _{TTQ} (\mu) \bigr\}.
  \end{multline*}
    It is easy to show that $\Sigma_l$ is a Lagrangian submanifold of $(T^*TQ, \omega_{TQ})$ (see \cite{Tu}).
Alternatively, we can write $\Sigma_l$ as
$$ \Sigma_l=\{\mu \in T ^\ast TQ\mid i_C^*\mu=dl \} = (\mathrm{d}L+\nu^{*}(C))\mid_{C}$$
with some abuse of notation. Here, $L: TQ\rightarrow \mathbb{R}$ is an arbitrary extension of $l$ to $TQ$ (that is $l\circ i_C=L$) and
$\nu^{*}(C)$  is the {\bf conormal bundle} of $C$:
\[
\nu^{*}(C)=\left\{\nu\in T^{*}TQ\big|_{C}\,\mid\,\langle \nu, T_{v}C\rangle=0 \hbox{  where } v=\pi_{TTQ}(\nu)\right\}.
\]
Therefore a curve $\mu: [t_0,t_1]\rightarrow \Sigma_l$ will be written
as $\mu(t)=\mathrm{d}L(\dot{\sigma}(t))+\nu(t),$ where $\nu(t)\in \left(\nu^{*}(C)\right)\mid_{\dot{\sigma}(t)}$ and $\dot{\sigma}(t)\in C\subseteq TQ,$
then
\begin{eqnarray*}
 \displaystyle 0 & =& \int_{t_0}^{t_1}\left\langle\mu(t),X^T(t, \dot{\sigma}(t))\right\rangle\ \mathrm{d}t=\int_{t_0}^{t_1}\left\langle\mathrm{d}L(\dot{\sigma}(t))+\nu(t),X^T(t, \dot{\sigma}(t))\right\rangle\ \mathrm{d}t\\
  &=&\int_{t_0}^{t_1}
  \left(\frac{d}{dt}\langle\mathrm{d}L + \nu,X^V\rangle(t, \ddot{\sigma}(t)) -\langle{\mathcal E}_{\mathrm{d}L+\nu} (\ddot{\sigma}(t)),
  X(t, \sigma(t))\rangle\right)\, \mathrm{d}t\\
  &=&-\int_{t_0}^{t_1}
  \langle\mathcal{E}_{\mathrm{d}L+\nu} (\ddot{\sigma}(t)), X(t,\sigma(t))\rangle\ \mathrm{d}t + \left[\langle\mathrm{d}L,X^V\rangle(t, \dot{\sigma}(t))+\langle \nu(t),
  X^V(t,\dot{\sigma}(t))\rangle\right]\Big |_{t_0}^{t_1}
  \end{eqnarray*}

Then, the equations of motion of the constrained variational problem are
\begin{eqnarray}
{\mathcal E}_{{\mathrm{d}L+\nu}}&=&0 \label{plm1}\\
\dot{\sigma}(t)&\in& C,\ \forall t\in[t_0,t_1], \label{plm2}
\end{eqnarray}
where a solution is a pair $(\sigma, \nu)$ with $\sigma: I\rightarrow Q$ and $\nu(t)\in \left[\nu^{*}(C)\right]_{\dot{\sigma}(t)}$.

Working locally,  assume that we have fixed local  constraints such that they determine $C$ by their vanishing, i.e., $\phi^{\alpha}(q,\dot q)=0$, $1\leq\alpha\leq m$, where
$m=\hbox{codim }  C$.
Therefore
\[
\left[\nu^{*}(C)\right]\mid_{\dot{\sigma}(t)}=\hbox{span }\{ \mathrm{d}\phi^{\alpha}(\dot{\sigma}(t))\}
\]
and, in consequence, $\nu(t)=\lambda_{\alpha}(t)\mathrm{d}\phi^{\alpha}(\dot{\sigma}(t))$ for some Lagrange multipliers $\lambda_{\alpha},$ to be determined. Then Equations (\ref{plm1}) and (\ref{plm2}) are now rewritten as
\begin{equation}\label{ecuacvakon}
\left\{
 \begin{array}{rcl}
{\mathcal E}_{\mathrm{d}L+\lambda_{\alpha} d\phi^{\alpha}} &=&0\\
\phi^{\alpha}(\dot{\sigma}(t))&=&0
\end{array}
\right.
\end{equation}
or, equivalently,
\begin{eqnarray*}
\displaystyle\frac{d}{dt}\left(\frac{\partial L}{\partial\dot{q}^i}+\lambda_{\alpha}\frac{\partial \phi^{\alpha}}{\partial\dot{q}^i} \right)-\frac{\partial L}{\partial q^i}-\lambda_{\alpha}\frac{\partial \phi^{\alpha}}{\partial {q}^i}&=&0\\
\phi^{\alpha}(q^i, \dot{q}^i)&=&0\; ,
\end{eqnarray*}
which are the equations of motion  for a constrained variational problem.

 Choosing an arbitrary frame $\{Y_i\}$ instead of the standard coordinate one, we im\-me\-dia\-te\-ly deduce that the equations of motion for the constrained variational problem are
\begin{eqnarray*}\label{gralvakonnoncanon}
\dot{q}^i-\rho^i_j y^j&=&0\\
\frac{d}{dt}\left(\frac{\partial (\tilde{L}+\lambda_{\alpha}\tilde{\phi}^{\alpha})}{\partial y^i}\right)
 -\rho^j_i\frac{\partial (\tilde{L}+\lambda_{\alpha}\tilde{\phi}^{\alpha})}{\partial q^j}+{\mathcal C}_{ij}^k
y^j\frac{\partial  (\tilde{L}+\lambda_{\alpha}\tilde{\phi}^{\alpha})}{\partial y^k}&=& 0\\
\tilde{\phi}^{\alpha}(q^i, y^i)&=&0,
\end{eqnarray*}
where $\tilde{\phi}^{\alpha}(q^i, y^i)=0$ are the constraint functions determining $C$ in terms of new coordinates $(q^i,y^i)$ and $\tilde{L}(q^i,y^i)=L(q^i,\rho^i_jy^j)$.

\vspace{.2cm}
An alternative way to describe the equations of motion in this case is related with the description
$\Sigma_l=\{\mu \in T ^\ast TQ\mid i_C^*\mu=\mathrm{d}l\}$, where we assume that the constraint functions are locally expressed as follows:
$\dot{q}^{\alpha}=\Phi^{\alpha}(q^i,\dot{q}^a)$, $1\leq\alpha\leq m$, $m+1\leq a\leq \dim Q.$ Hence, $i_C:C\hookrightarrow TQ$ is written as
$i_C(q^i,\dot{q}^a)=(q^i,\dot{q}^a,\Phi^{\alpha}(q^i,\dot{q}^a))$ and, if we take an arbitrary 1-form $$\mu= \mu_i\mathrm{d}q^i+\tilde{\mu}_a\mathrm{d}\dot{q}^a+\tilde{\mu}_{\alpha}
   \mathrm{d}\dot{q}^{\alpha}$$ on $TQ,$
then
   $$\displaystyle i^*_C\mu=\mu_i\mathrm{d}q^i+\tilde{\mu}_a\mathrm{d}\dot{q}^a+
   \tilde{\mu}_{\alpha}\frac{\partial\Phi^{\alpha}}{\partial q^i}\mathrm{d}q^i+\tilde{\mu}_{\alpha}
   \frac{\partial\Phi^{\alpha}}{\partial\dot{q}^a}\mathrm{d}\dot{q}^a.$$

Since  $i_C^*\mu=\mathrm{d}l,$ then
\begin{eqnarray*}
\mu_i&=&\frac{\partial l}{\partial q^i}-\tilde{\mu}_{\alpha}\frac{\partial\Phi^{\alpha}}{\partial q^i}\\
\tilde{\mu}_a&=&\frac{\partial l}{\partial\dot{q}^a}-\tilde{\mu}_{\alpha}
   \frac{\partial\Phi^{\alpha}}{\partial\dot{q}^a}\; .
\end{eqnarray*}
Observe that we are naturally describing $\Sigma_l$ with coordinates $(q^i, \dot{q}^a, \tilde{\mu}_{\alpha})$. Thus, applying the generalized variational
calculus to $\Sigma_l,$ we arrive to an alternative but equivalent description of the constrained variational calculus by the equation
   \begin{eqnarray*}
\displaystyle 0 & =& \int_{t_0}^{t_1}\left\langle\mu(t),X^T(t, \dot{\sigma}(t))\right\rangle\ \mathrm{d}t\\
   &=& \displaystyle\int_{t_0}^{t_1}\left\langle\left(\frac{\partial l}{\partial q^i}-\tilde{\mu}_{\alpha}\frac{\partial\Phi^{\alpha}}{\partial q^i}\right)\mathrm{d}q^i+\left(\frac{\partial l}{\partial \dot{q}^a}-\tilde{\mu}_{\alpha}\frac{\partial\Phi^{\alpha}}{\partial \dot{q}^a}\right)\mathrm{d}\dot{q}^a+
\tilde{\mu}^{\alpha}\mathrm{d}\dot{q}^{\alpha}, X^T(t, \dot{\sigma}(t))\right\rangle \mathrm{d}t
\end{eqnarray*}
from which we easily derive the equations
\begin{eqnarray*}
\frac{d}{dt}\left(\frac{\partial l}{\partial \dot{q}^a}-\tilde{\mu}_{\alpha}\frac{\partial\Phi^{\alpha}}{\partial \dot{q}^a}\right)&=&\frac{\partial l}{\partial q^a}-\tilde{\mu}_{\alpha}\frac{\partial\Phi^{\alpha}}{\partial q^a}\\
\frac{d\tilde{\mu}_{\alpha}}{dt}&=&\frac{\partial l}{\partial q^{\alpha}}-\tilde{\mu}_{\alpha}\frac{\partial\Phi^{\alpha}}{\partial q^{\alpha}}\\
\frac{d\dot{q}^{\alpha}}{dt}&=&\Phi^{\alpha}(q^i, \dot{q}^a).
\end{eqnarray*}

These equations are obtained in \cite{CdDdLM03} using variational techniques and introducing an ansatz in the deduction that now is clarified in the context of the generalized variational calculus.

\vspace{.3cm}
In  coordinates $(q^i,y^i),$ assuming that the constraint submanifold $C$ is locally given by the vanishing of the constraints $y^{\alpha}=\Phi^{\alpha}(q^i,y^a),$ we have $i_C:C\hookrightarrow TQ$ given by $i_C(q^i,y^a)=(q^i,y^a,\Phi^{\alpha}(q^i,y^a))$ and we take \newline $\Sigma=\{\mu=\mu_i\mathrm{d}q^i+\tilde{\mu}_i\mathrm{d}y^i=\mu_i\mathrm{d}q^i+\tilde{\mu}_a\mathrm{d}y^a+\tilde{\mu}_{\alpha}\mathrm{d}y^{\alpha}
:i^*_C\mu=\mathrm{d}l\},$ where $y^{\alpha}=\Phi^{\alpha}(q^i,y^a).$ Thus,
$$\mu=\left(\frac{\partial l}{\partial q^i}-\tilde{\mu}_{\alpha}\frac{\partial\Phi^{\alpha}}{\partial q^i}\right)\mathrm{d}q^i+\left(\frac{\partial l}{\partial y^a}-\tilde{\mu}_{\alpha}\frac{\partial\Phi^{\alpha}}{\partial y^a}\right)\mathrm{d}y^a+\tilde{\mu}_{\alpha}\mathrm{d}y^{\alpha}.$$

From (\ref{esubmunoncan}) we have that the equations of the generalized variational calculus in this case are

\begin{eqnarray}\label{ecuacVakonCuasi}
\displaystyle\frac{d}{dt}\left(
\frac{\partial l}{\partial \dot{y}^a}-\tilde{\mu}_{\alpha}\frac{\partial\Phi^{\alpha}}{\partial \dot{y}^a}\right)-\left(\frac{\partial l}{\partial q^i}-\tilde{\mu}_{\alpha}\frac{\partial\Phi^{\alpha}}{\partial q^i}\right)\rho_a^i+\tilde{\mu}_k\mathcal{C}^k_{aj}y^j &=& 0\\
\displaystyle\frac{d}{dt}\tilde{\mu}_{\alpha}-\left(\frac{\partial l}{\partial q^i}-\tilde{\mu}_{\alpha}\frac{\partial\Phi^{\alpha}}{\partial q^i}\right)\rho_{\alpha}^i+\tilde{\mu}_k\mathcal{C}^k_{\alpha j}y^j &=& 0\\
\displaystyle y^{\alpha}&=&\Phi^{\alpha}(q^i, y^a),
\end{eqnarray}
for $1\leq i\leq n,$ $1\leq\alpha\leq m$ and $1\leq a\leq n-m,$ where $\dim C=n-m.$ Then, using the expression for $\tilde{\mu}_i$ and $y^j,$ we obtain the
following system of equations for vakonomic mechanics
\begin{eqnarray*}\label{ecuacvakonnoncanon}
0&=&\displaystyle\frac{d}{dt}\left(
\frac{\partial l}{\partial y^a}-\tilde{\mu}_{\alpha}\frac{\partial\Phi^{\alpha}}{\partial y^a}\right)-\left(\frac{\partial l}{\partial q^i}-\tilde{\mu}_{\alpha}\frac{\partial\Phi^{\alpha}}{\partial q^i}\right)\rho_a^i+
\left(\frac{\partial l}{\partial y^c}-\tilde{\mu}_{\alpha}\frac{\partial\Phi^{\alpha}}{\partial y^c}\right)\mathcal{C}^c_{a\beta}\Phi^{\beta} \\
&&+ \left(\frac{\partial l}{\partial y^c}-\tilde{\mu}_{\alpha}\frac{\partial\Phi^{\alpha}}{\partial y^c}\right)\mathcal{C}^c_{a b}y^b +\tilde{\mu}_{\gamma}\mathcal{C}^{\gamma}_{a b}y^b +\tilde{\mu}_{\gamma}\mathcal{C}^{\gamma}_{a \beta}\Phi^{\beta} \\
0&=&\displaystyle\frac{d}{dt}\tilde{\mu}_{\alpha}-\left(\frac{\partial l}{\partial q^i}-\tilde{\mu}_{\delta}\frac{\partial\Phi^{\delta}}{\partial q^i}\right)\rho_{\alpha}^i+\tilde{\mu}_{\gamma}\mathcal{C}^{\gamma}_{\alpha b}y^b +\tilde{\mu}_{\gamma}\mathcal{C}^{\gamma}_{\alpha \beta}\Phi^{\beta} \\
& &
+\left(\frac{\partial l}{\partial y^c}-\tilde{\mu}_{\delta}\frac{\partial\Phi^{\delta}}{\partial y^c}\right)\mathcal{C}^{c}_{\alpha b}y^b+
 \left(\frac{\partial l}{\partial y^c}-\tilde{\mu}_{\delta}\frac{\partial\Phi^{\delta}}{\partial y^c}\right)\mathcal{C}^{c}_{\alpha \beta}\Phi^{\beta} \\
\displaystyle \dot{q}^i&=&\rho_a^i y^a+\rho^i_{\alpha}\Phi^{\alpha}.
\end{eqnarray*}

These equations coincide with the ones derived in \cite{DianaRestringido}.

\subsubsection{Sub-Riemannian geometry.}

Sub-Riemannian geometry is a  generalization of Riemannian geometric where the Riemannian metric is only defined on a vector subbundle of the tangent bundle
to the manifold, instead on the full manifold. The notion of length is only assigned to a particular  subclass of curves, that is,  curves with tangent vectors belonging to the vector subbundle for each point.
 More precisely, we consider a  manifold $Q$ equipped with a smooth distribution ${\mathcal D}$ of constant rank. A sub-Riemannian metric on ${\mathcal D}$
 consists of a positive definite quadratic form $g_q$  on ${\mathcal D}_q$ smoothly varying in $q\in Q$. We will say that a piecewise smooth curve $\sigma: [t_0, t_1]\rightarrow Q$ is admissible if
 $\dot{\sigma}(t)\in {\mathcal D}_{\sigma(t)}$, for all $t\in [t_0, t_1]$. We define its  length as follows
 \[
 \hbox{lenght}(\sigma)=\int_{t_0}^{t_1}\sqrt{g(\dot{\sigma}(t),\dot{\sigma}(t)}\, \mathrm{d}t.
 \]
  From this definition, we have a notion of distance between two points  $x,y\in Q$ as
$\hbox{dist } (x, y)=\hbox{inf}_{\sigma}\hbox{lenght}(\sigma)$.
It is finite if there  exists admissible curves $\sigma$ connecting $x$ and $y$; in another case, the distance is considered infinite.
  A curve which realizes the distance between two points is called a \emph{minimizing sub-Riemannian geodesic}. It is clear that the problem of finding
   minimizing sub-Riemannian geodesics is
exactly the same as the vakonomic problem determined by the restricted Lagrangian $l: {\mathcal D}\rightarrow {\mathbb R}$ defined by
   $l(v_q)=\frac{1}{2}g_q(v_q, v_q),$ where $v_q\in {\mathcal D}_q$.

  Now, we will see a particular example of  sub-Riemannian geometry.
We consider a local sub-Riemannian problem given by $(U, {\mathcal D},g),$ where $U$ is an open set in $\mathbb{R}^3$ con\-tai\-ning $(0, 0, 0),$
 ${\mathcal D}$ is the distribution $\ker(\omega)$ being $\omega$ the Martinet 1-form
  $\displaystyle \mathrm{d}q^3-\frac{(q^2)^2}{2}\mathrm{d}q^1$ and $q=(q^1,q^2,q^3)$ are the coordinates. The sub-Riemannian metric $g$
 is defined on ${\mathcal D}$ by $a(q)\mathrm{d}(q^1)^2+2b(q)\mathrm{d}q^1\mathrm{d}q^2+c(q)\mathrm{d}(q^2)^2$
 but, for simplicity,  we assume that $a(q)=1,b(q)=0$ and $c(q)=1/2.$
So, in our notation, \newline $C=\displaystyle\left\{(q^1,q^2,q^3,\dot{q}^1,\dot{q}^2,\dot{q}^3)\in TU\equiv U\times {\mathbb R}^3\; |\;  \dot{q}^3=\frac{(q^2)^2}{2}\dot{q}^1\right\}$.
Therefore, we have  coordinates $(q^1,q^2,q^3,\dot{q}^1,\dot{q}^2)$ in $C$ with inclusion
 $\displaystyle i_C:(q^1,q^2,q^3,\dot{q}^1,\dot{q}^2)\mapsto\left(q^1,q^2,q^3,\dot{q}^1,\dot{q}^2,\frac{(q^2)^2}{2}\dot{q}^1\right)$.
 The  Lagrangian function $l:C\rightarrow\mathbb{R}$ is given by
 $\displaystyle l(q^1,q^2,q^3,\dot{q}^1,\dot{q}^2)=\frac{1}{2}\left((\dot{q}^1)^2+(\dot{q}^2)^2\right),$ $\Sigma=\{\mu\in T^*TU: i^*\mu=\mathrm{d}l\},$ and a  1-form $\mu\in\Sigma\subseteq T^*TU$ can be written as \newline $\mu= \displaystyle -q^2\dot{q}^1\tilde{\mu}_3\mathrm{d}q^2+\left(2\dot{q}^1-\frac{(q^2)^2}{2}\tilde{\mu}_3\right)\mathrm{d}\dot{q}^1
+2\dot{q}^2\mathrm{d}\dot{q}^2+\tilde{\mu}_3\mathrm{d}\dot{q}^3.$

If we consider the adapted basis $\displaystyle\{Y_1,Y_2,Y_3\}=\left\{\frac{\partial}{\partial q^2},\frac{\partial}{\partial q^1}+\frac{(q^2)^2}{2}\frac{\partial}{\partial q^3},\frac{\partial}{\partial q^3}\right\}$
of vector fields on $Q$, we induce coordinates  $\{y^1,y^2,y^3\}$ where now $C$ is determined by the constraint $y^3=0.$ In this case, we obtain that
$\mathcal{C}^3_{12}=y=-\mathcal{C}^3_{21},$ $\rho_1^2=\rho_2^1=\rho_3^3=1$ and $\rho_2^3=y^2/2.$ Hence, applying equations  (\ref{ecuacVakonCuasi}), we obtain
\begin{equation*}
\left\{
\begin{array}{rl}
\displaystyle \dot{y}^1+\tilde{\mu}_3q^2y^2&=0 \\
\dot{y}^2-\tilde{\mu}_3q^2y^1 &=0 \\
\dot{\tilde{\mu}}_3 &=0 \\
\end{array}
\right.\qquad\text{and}\qquad
\left.
\begin{array}{rl}
\dot{q}^1 &=y^2 \\
\dot{q}^2 &=y^1 \\
\dot{q}^3 &=\displaystyle\frac{(q^2)^2}{2}y^2
\end{array}.
\right.
\end{equation*}
These equations coincide with the ones obtained in \cite{BonnardChyba}.

\begin{remark}
{\rm
It is interesting to note that our formalism is also adapted to the study of abnormal solutions of sub-Riemannian geometry (see \cite{MR1867362}). For a complete study of regular and normal solutions, it is only necessary
to consider the subset $\Sigma=\Sigma_l\cup \nu^{*}({\mathcal D})$.
}
\end{remark}

\subsubsection{Higher-order Lagrangian systems}
In the case in which we have a higher-order Lagrangian $L:T^{(k)}Q\rightarrow\mathbb{R}$, that is, a Lagrangian depending on higher-order derivatives (positions, velocities, acelerations and so on), we can also apply the generalized variational calculus. As in Section \ref{sect:fibtangente},
we know that we can see $T^{(k)}Q$ as a submanifold of $TT^{(k-1)}Q$, using the inclusion $j_k:T^{(k)}Q\hookrightarrow TT^{(k-1)}Q$ (see \cite{ManuelCampos}). With this point of view, we can see any higher-lagrangian problem as a constrained variational problem where we take
$\Sigma_L=\{\mu\in T^*TT^{(k-1)}Q\; \mid \; j_k^*\mu=\mathrm{d}L\}.$

In this case, a curve $\sigma:I\rightarrow Q$ is a solution of the higher-order variational problem determined by $L:T^{(k)}Q\rightarrow\mathbb{R}$ if its lift $\sigma^{(k-1)}:I\rightarrow T^{(k-1)}Q$ is solution of the generalized variational problem determined by $\Sigma_L$.
Schematically, we have
\begin{equation*}
\xymatrix{
&\Sigma\ar[d]_{(\pi_{TTQ})_{|\Sigma}}\ar@{^{(}->}[rr]&&T^*TT^{(k-1)}Q\ar[d]\\
&T^{(k)}Q\ar@{^{(}->}[rr]^{j_{k-1}}&&TT^{(k-1)}Q\ar[d]_{\tau_{T^{(k-1)}Q}}\\
\mathbb{R}\ar[rrr]_{\sigma^{(k-1)}}\ar[ur]^{\dot{\sigma}^{(k-1)}}\ar@/^3.5pc/[uur]^{\mu} & &&T^{(k-1)}Q\\
}
\end{equation*}

\subsection{Constrained Hamiltonian mechanics}
We can combine the techniques developed in sections \ref{mecHam} and \ref{calcvariacrestringido} to produce the equations for contrained Hamiltonian mechanics (topic treated in \cite{dirac1,dirac2,dirac3,MR1692290}). For that, consider a submanifold $\mathcal{M}$ of $T^*Q$ with canonical inclusion $i_{\mathcal{M}}: {\mathcal M}\hookrightarrow T^*Q,$ and a function $h:\mathcal{M}\rightarrow\mathbb{R}.$

  Let's take the (Lagrangian-)submanifold $S_h=\{\mu\in T^*T^*Q\; \mid\; i_{\mathcal{M}}^*\mu=\mathrm{d}h\}$ of $T^*T^*Q.$ Using the antisymplectomorphism $\mathcal{R}:T^*T^*Q\rightarrow T^*TQ$ (see Section \ref{mecHam}), we induce a submanifold $\Sigma_h=\mathcal{R}(S_h)$ of $T^*TQ$ and, therefore, we can use the generalized variational calculus to characterize the solutions of the constrained Hamiltonian problem.

In local coordinates, if $\mathcal{M}$ is determined by the vanishing of constraints $\Phi^{\alpha}(q,p)=0$ and  $H:T^*Q\rightarrow\mathbb{R}$ is an arbitrary extension of $h:\mathcal{M}\rightarrow\mathbb{R},$ (i.e., $h\circ i_{\mathcal{M}}=H\mid_{\mathcal M}$), then $$S_h=\{\mu\in T^*T^*Q\; \mid\; \mu=\mathrm{d}H+\lambda_{\alpha}\mathrm{d}\Phi^{\alpha},\; \Phi^{\alpha}=0\}.$$
If we choose local coordinates $\displaystyle (q^i,p_i,\eta_i,\tilde{\eta}^i)$  in $T^*T^*Q,$ then $$\displaystyle S_h=\left\{(q^i,p_i,\eta_i,\tilde\eta^i)\; \mid\; \eta_i=\frac{\partial H}{\partial q^i}+\lambda_{\alpha}\frac{\partial\Phi^{\alpha}}{\partial q^i},\ \tilde\eta^i=\frac{\partial H}{\partial p_i}+\lambda_{\alpha}\frac{\partial\Phi^{\alpha}}{\partial p_i}\right\}$$ and $\Sigma_h \subseteq T^*TQ$ is given by
$$\Sigma_h=\left\{(q^i,\dot{q}^i,\mu_i,\tilde\mu_i)\; \mid\;  \dot{q}^i=\frac{\partial H}{\partial p_i}+\lambda_ {\alpha}\frac{\partial\Phi^{\alpha}}{\partial p_i},\mu_i=-\frac{\partial H}{\partial q^i}-\lambda_{\alpha}\frac{\partial\Phi^{\alpha}}{\partial q^i}\right\}.$$

A curve $\mu:I\rightarrow\Sigma$ is a solution of the generalized variational problem induced by $h:\mathcal{M}\rightarrow\mathbb{R}$ if $\mu(t)=(q^i(t),\dot{q}^i(t),\dot{\mu}_i(t),\tilde{\mu}_i(t))$ verifies the equations

\begin{equation}\label{E-L}
\left\{
\begin{array}{rl}
\mathcal{E}_{\mu}&=\displaystyle\frac{\mathrm{d}\tilde\mu_i(t)}{\mathrm{d}t}+\frac{\partial H}{\partial q^i}(q(t),\tilde\mu(t))+\lambda_{\alpha}(t)\frac{\partial\Phi^{\alpha}}{\partial q^i}(q(t), \tilde\mu(t))=0\\
\dot{q}^i(t)&=\displaystyle\frac{\partial H}{\partial p_i}(q(t),\tilde\mu(t))+\lambda_{\alpha}(t)\frac{\partial\Phi^{\alpha}}{\partial p_i}(q(t),\tilde\mu(t)) \\
\displaystyle \Phi^{\alpha}(q^i(t),\tilde\mu(t))&=0.
\end{array}
\right.
\end{equation}
\vspace{.2cm}
Taking the time derivative of the last equation, we obtain

$$\frac{\partial\Phi^{\alpha}}{\partial q^i}(q(t),\tilde\mu(t))\dot{q}^i(t)+\frac{\partial\Phi^{\alpha}}{\partial p_i}(q(t),\tilde\mu(t))\dot{\tilde\mu}_i(t)=0$$

and, using equations (\ref{E-L}), we get a new constraint equation:
\begin{equation}\label{nuevaconstraint}
 \frac{\partial\Phi^{\alpha}}{\partial q^i}\left(\frac{\partial H}{\partial p_i}+\lambda_ {\alpha}\frac{\partial\Phi^{\alpha}}{\partial p_i}\right)- \frac{\partial\Phi^{\alpha}}{\partial p_i}\left(\frac{\partial H}{\partial q^i}+\lambda_ {\alpha}\frac{\partial\Phi^{\alpha}}{\partial q^i}\right)=0.
\end{equation}
Proceeding further, we would derive the classical Dirac-Bergmann constraint algorithm (see \cite{bergmann,dirac3}).

\subsection{Optimal Control Theory}

Generally speaking, an optimal control problem from the differential geometric viewpoint is given by a vector field depending on some parameters called controls, some boundary conditions and a cost function whose integral must be either minimized or maximized.
Concretely,  an \emph{optimal control problem} $(U, Q,\Gamma, { L})$ is given by a control bundle
$\tau_{U,Q}\colon {U} \rightarrow Q$, a vector field $\Gamma$ along the control bundle
projection $\tau_{{U},Q}$, a cost function ${ L} \colon { U} \rightarrow \mathbb{R}$ whose functional
must be minimized, and some endpoint or boundary conditions that must be satisfied at initial and/or final time.
By definition, the vector field  $\Gamma$ along $\tau_{ U,Q}$ verifies $\tau_{TQ}\circ \Gamma=\tau_{{ U},Q}$. We have the diagram
\begin{equation*}\xymatrix{&& TQ \ar[d]^{\txt{\small{$\tau_{TQ}$}}} \\ {U}
\ar[rr]^{\txt{\small{$\tau_{{U},Q}$}}} \ar[rru]^{\txt{\small{$\Gamma$}}}& & Q}
\end{equation*}

Locally, $\Gamma(q^i, u^a)=(q^i, \Gamma^i(q, u))$ which defines the control equations $\dot{q}^i=\Gamma^i(q,u)$.
From the optimal control data  $( U, Q, \Gamma, { L}),$ we construct Pontryagin's hamiltonian \newline $H:T^*Q\times_Q  U\longrightarrow \mathbb{R}$ given by
\begin{equation}\label{pont}
H(\alpha_q, u_q)=\langle \alpha_q, \Gamma(u_q)\rangle-L(u_q),
\end{equation}
where $u_q\in U_q$ and $\alpha_q\in T_q^*Q$.
In coordinates, $H(q^i, p_i, u^a)=p_i\Gamma^i(q, u)-L(q, u)$.
The usual technique to solve an optimal control problem is Pontryagin's Maximum Principle (see, for instance, \cite{MR2551483,MR0186436}), which
provides us with a set of necessary conditions for optimality.

 The optimal control solutions can be also characterized using the generalized variational calculus. For that, we define the subset of $T^*TQ$:
\[
\Sigma=\{\mu\in T^*TQ\; \mid\; \Gamma^*\mu=dL\}\, .
\]
Observe that it is not in general a submanifold of $T^*TQ$.

 The solutions are determined by the conditions:
$\mathcal{E}_{\mu}=0,$ with $\mu: I\rightarrow \Sigma,$ and $\pi_{TTQ}(\mu(t))=\dot{\sigma}(t)\in \Gamma(U)=C.$

Locally, if we take an arbitrary element $\mu=\mu_i\mathrm{d}q^i+\tilde\mu_j\mathrm{d}\dot{q}^j\in T^*TQ,$
then $\mu\in \Sigma$ if:
\begin{eqnarray*}
\Gamma^*\mu&=& \Gamma^*\left(\mu_i\mathrm{d}q^i+\tilde\mu_j\mathrm{d}\dot{q}^j\right)\\
&=&\mu_i\mathrm{d}q^i+\tilde\mu_j\mathrm{d}(\Gamma^j(q,u))\\
&=&\left(\mu_i+\frac{\partial \Gamma^j}{\partial q^i}\tilde\mu_j\right)\mathrm{d}q^i+\tilde\mu_j\frac{\partial \Gamma^j}{\partial u^{a}}\mathrm{d}u^{a}\\
&=&\frac{\partial L}{\partial q^i}\mathrm{d}q^i+\frac{\partial L}{\partial u^{a}}\mathrm{d}u^{a}\\
&=&dL
\end{eqnarray*}

From Definition \ref{solution-m}, we have that a  curve $\mu: I \rightarrow \Sigma$ is a solution of the generalized variation problem determined by $\Sigma$ if it verifies the following set of equations
\begin{equation}
   \left\{
   \begin{array}{rcl}
    \displaystyle \frac{d}{dt}\tilde\mu_i(t)-\mu_i(t)=0\; \\
      \displaystyle \mu_i(t)+\tilde\mu_j\frac{\partial \Gamma^j}{\partial q^i}(q(t), u(t))\tilde\mu_j(t)=\frac{\partial L}{\partial q^i}(q(t), u(t))\\
      \displaystyle \tilde\mu_j(t)\frac{\partial 	\Gamma^j}{\partial u^{a}}(q(t), u(t))=\frac{\partial L}{\partial u^{a}}(q(t), u(t))\\
      \displaystyle \dot{q}^i(t)=\Gamma^i(q(t),u(t))\; .
   \end{array}
   \right.
\end{equation}

Replacing the expression of $\mu_i$ from the second equation in the first one, we obtain the following system:

\begin{equation}\label{H-PPrinciple}
   \left\{
   \begin{array}{rcl}
      \displaystyle \frac{d \tilde\mu_i}{dt}(t)=\frac{\partial L}{\partial q^i}(q(t), u(t))-\tilde\mu_j(t)\frac{\partial \Gamma^j}{\partial q^i}(q(t), u(t))\tilde\mu_j(t)\\
      \displaystyle \frac{\partial L}{\partial u^{a}}(q(t), u(t))-\tilde\mu_j(t)\frac{\partial \Gamma^j}{\partial u^{a}}(q(t), u(t))=0\\
      \displaystyle \frac{d{q}^i}{dt}(t)=\Gamma^i(q(t), u(t)),
   \end{array}
   \right.
\end{equation}
that in terms of Pontryagin's hamiltonian  $H:T^*Q\times_Q  U\longrightarrow \mathbb{R}$ defined in (\ref{pont}), can be written as

\begin{equation}\label{H-PPrinciple}
   \left\{
   \begin{array}{rcl}
      \displaystyle \frac{d\tilde\mu_i}{dt}(t)=-\frac{\partial H}{\partial q^i}(q(t), u(t))\\
      \displaystyle \frac{\partial H}{\partial u^{a}}(q(t), u(t))=0\\
      \displaystyle \frac{dq^i}{dt}(t)=\frac{\partial H}{\partial \tilde\mu_i}(q(t), u(t))
      \end{array}
   \right.
\end{equation}
to regain Hamilton-Pontryaguin's conditions of extremality.

In a coordinate system $\{(q^i,y^j)\}$ adapted to an arbitrary frame, we have $\mu=\mu_i\mathrm{d}q^i+\tilde{\mu}_j\mathrm{d}y^j$ and $\Gamma(q^i,u^a)=(q^i,\Gamma^j(q^i,u^a)).$ Then we obtain that the conditions for $\mu$ belonging to $\Sigma$ are
\begin{equation}
   \left\{
   \begin{array}{rcl}
      \displaystyle \mu_i &=&\displaystyle\frac{\partial L}{\partial q^i}-\tilde{\mu}_j\frac{\partial\Gamma^j}{\partial q^i}\\
      \displaystyle \tilde{\mu}_j\frac{\partial \Gamma^j}{\partial u^a}& =&\displaystyle\frac{\partial L}{\partial u^a}\\
      \displaystyle y^i &=&\Gamma^i(q^i, u^a),
      \end{array}
   \right.
\end{equation}
and using expression  (\ref{esubmunoncan}), we deduce the following system of equations
\begin{equation}
   \left\{
   \begin{array}{rcl}
      \displaystyle \frac{d}{dt}\tilde{\mu}_i &=&\displaystyle\left(\frac{\partial L}{\partial q^j}-\tilde{\mu}_k\frac{\partial \Gamma^k}{\partial q^j}\right)\rho_i^j-\mathcal{C}^k_{ij}y^j\tilde{\mu}_k \\
      \displaystyle \tilde{\mu}_i\frac{\partial\Gamma^i}{\partial u^a}& =&\displaystyle \frac{\partial L}{\partial u^a} \\
       \dot{q}^i &=&\displaystyle \rho^i_jy^j=\rho_j^i\Gamma^j,
   \end{array}
   \right.
\end{equation}
where the last equation is the condition for admissibility.

Hence, in terms of Pontryagin's Hamiltonian, the equations of optimal control obtained by using the  generalized variational calculus are
\begin{equation}
   \left\{
   \begin{array}{rcl}
      \displaystyle \frac{{d}}{{d}t}\tilde{\mu}_i &=&\displaystyle -\frac{\partial H}{\partial q^j}\rho_i^j-\mathcal{C}^k_{ij}\Gamma^j\tilde{\mu}_k \\
      \displaystyle \frac{\partial H}{\partial u^a}& =&0 \\
       \dot{q}^i &=&\displaystyle \rho^i_j\frac{\partial H}{\partial\tilde{\mu}_j}
      \end{array}
   \right.
\end{equation}

where $H(qî,\tilde{\mu}_j,u^a)=\tilde{\mu}_j\Gamma^j(q,u)-L(q,u).$

\section{\textbf{\textsc{Generalized variational calculus on skew-symmetric Lie algebroids}}} \label{sect:algebroides}
Now, we will show an extension of the generalized variational calculus to other different system of great interest in mechanics: reduced systems and nonholonomic systems.
In many cases,  Lagrangian  or Hamiltonian systems admit a group of symmetries and it is possible to reduce the original system to a new one defined on a reduced space with less degrees of freedom or, in other case, the phase space is reduced due to the presence of nonholonomic constraints.
The theory of Lie algebroids or, more generally, skew-symmetric algebroids, provides an unifying framework for all these systems (see \cite{dLMMalgebroides, GraGra08,martinez01,weinstein78}).

\subsection{Skew-symmetric algebroids} A \textbf{skew-symmetric algebroid structure} on a vector bundle $\tau_D:D\rightarrow Q$ is a $\mathbb{R}$-bilinear bracket $\lcf .,.\rcf_D:\Gamma(\tau_D)\times\Gamma(\tau_D)\rightarrow\Gamma(\tau_D)$ on the space $\Gamma(\tau_D)$ of sections of $\tau_D$ and a vector bundle morphism $\rho_D:D\rightarrow TQ,$ which is called $\textbf{anchor map},$ such that:
\begin{enumerate}
   \item $\lcf .,.\rcf_D$ is skew-symmetric; that is, $\lcf X,Y\rcf_D=-\lcf Y,X\rcf_D,$ for $X,Y\in\Gamma(\tau_D).$
   \item If we denote by $\rho_D:\Gamma(\tau_D)\rightarrow\mathfrak{X}(Q)$ the morphism of $C^{\infty}(Q)$-modules induced by the anchor map, then
   $\lcf X,fY\rcf_D=f\lcf X,Y\rcf_D+\rho_D(X)(f)Y,$ for $X,Y\in\Gamma(\tau_D)$ and $f\in C^{\infty}(Q).$
\end{enumerate}

If the bracket $\lcf .,.\rcf_D$ satisfies the Jacobi identity, we have that the pair $(\lcf .,.\rcf_D,\rho_D)$ is a \textbf{Lie algebroid structure} on the vector bundle $\tau_D:D\rightarrow Q.$
Therefore, a Lie algebroid over a manifold $Q$ may be thought of as a
``generalized tangent bundle'' to $Q$. We will see some interesting examples where this structure appears. For more details see \cite{mackenzie87}.

\begin{itemize}
\item A \emph{finite dimensional real Lie algebra} $\mathfrak{g}$ where
$Q=\{q\}$ be a unique point. Then, we consider the vector bundle
$\tau_{\mathfrak{g}}:\mathfrak{g}\rightarrow \{ q\}.$ The sections of this
bundle can be identified with the elements of $\mathfrak{g}$ and,
therefore, we can consider the structure of
the Lie algebra $[\cdot,\cdot]_{\mathfrak{g}}$ as the Lie bracket. The anchor map
is $\rho\equiv 0.$ Then,
$(\mathfrak{g},[\cdot,\cdot]_{\mathfrak{g}},0)$ is a Lie
algebroid over a point.

\item A \emph{tangent bundle} of a manifold $Q$ (see Section \ref{sect:fibtangente}). The sections of the bundle $\tau_{TQ}:TQ\rightarrow Q$ are identified with the vector fields on
$Q,$ the anchor map $\rho:TQ\rightarrow TQ$ is the identity function and
the Lie bracket defined on $\Gamma(\tau_{TQ})$ is induced by the standard
Lie bracket of vector fields on $Q.$

\item Let $\phi:Q\times G\rightarrow Q$ be a right action of $G$ on the manifold $Q$ where
$G$  is a Lie group. The induced anti-homomorphism between the
Lie algebras $\mathfrak{g}$ and $\mathfrak{X}(Q)$ is given by
$\Delta:\mathfrak{g}\rightarrow\mathfrak{X}(Q),$ where
$\Delta(\xi)=\xi_Q$ is the infinitesimal generator of the
action for $\xi\in\mathfrak{g}.$
\par The vector bundle $\tau_{Q\times\mathfrak{g}}:Q\times\mathfrak{g}\rightarrow
Q$ is a Lie algebroid over $Q$ where the anchor map $\rho:Q\times\mathfrak{g}\rightarrow TQ$ is defined
as $\rho(q,\xi)=-\xi_{Q}(q),$ and
the bracket is given by the Lie algebra structure on
$\Gamma(\tau_{Q\times\mathfrak{g}})$ by
$$\lcf\hat{\xi},\hat{\eta}\rcf(q)=(q,[\xi,\eta])=\widehat{[\xi,\eta]},\ \hbox{ for } q\in Q,$$
where $\hat{\xi}(q)=(q,\xi),$ and $\hat{\eta}(q)=(q,\eta),$ with $\xi, \eta\in {\mathfrak g}$. This Lie
algebroid is know as \emph{Action Lie algebroid}.

\item Let $G$ be a Lie group and we assume that $G$ acts freely and properly on $Q$. We
denote by $\pi:Q\rightarrow \widehat{Q}=Q/G$ the associated principal
bundle. The tangent lift of the action gives a free and proper
action of $G$ on $TQ$ and we denote by $\widehat{TQ}=TQ/G$ the corresponding quotient
manifold. The quotient vector bundle
$\tau_{\widehat{TQ}}:\widehat{TQ}\rightarrow Q/G,$ where
$\tau_{\widehat{TQ}}([v_q])=\pi(q),$ is a Lie algebroid over
$Q/G.$ The Lie bracket is defined on $\Gamma(\tau_{TQ/G})$ and it is isomorphic to the Lie subalgebra of $G$-invariant vector
fields. Thus, the Lie bracket on $\widehat{TQ}$ is just the
bracket of $G-$invariant vector fields. The anchor map
$\rho:TQ/G\rightarrow T(Q/G)$ is given by $\rho([v_q])=T_{q}\pi(v_q).$
This Lie algebroid is called \emph{Atiyah algebroid}
associated with the principal bundle $\pi:Q\rightarrow\widehat{Q}.$
\end{itemize}

\vspace{.3cm}
Suppose that $(q^i)$ are local coordinates on $Q$ and that
$\{e_{A}\}$ is a local basis of the space of sections $\Gamma(\tau_D)$,
then
\begin{equation*} \lcf e_{A}, e_B\rcf_D =
{\mathcal C}^C_{AB} e_{C}, \ \ \ \rho_D (e_{A})=\rho_{A}^i \frac{\partial}{\partial q^i}  .\label{coeff-estruct-A}
\end{equation*}
The functions ${\mathcal C}^C_{AB},\ \rho_A^i\in C^{\infty}(Q)$
are called the \emph{local structure functions} of the
skew-symmetric algebroid on $\tau_D: D \rightarrow Q$.

A $\rho_D$\emph{-admissible curve} is a curve $\gamma: I\subseteq
\mathbb{R} \longrightarrow D$  such that
\[
\frac{{d}(\tau_D\circ \gamma)}{{ d}t}(t)=\rho_D(\gamma(t))\;
.
\]
 Let's define the set \begin{equation}\label{D2}
D^{(2)}=\{v\in TD\; |\; \rho_d(\tau_D(v))=T\tau_D(v)\}
\end{equation}
which will play a similar role to $T^{(2)}Q$ in Section \ref{sect:fibtangente}.
We can define $D^{(2)}$ in an alternative way. Considering two admissible curves $\gamma_1:I\rightarrow D$ and $\gamma_2:I\rightarrow D$ such that $\gamma_1(0)=\gamma_2(0),$ we say that $\gamma_1$ and $\gamma_2$ define the same equivalence class if and only if $\displaystyle \frac{\mathrm{d}\gamma_1}{\mathrm{d}t}(0)=\frac{\mathrm{d}\gamma_2}{\mathrm{d}t}(0).$ The set of these equivalence classes is precisely the set $D^{(2)}$ defined as in (\ref{D2}). We will denote by $[\gamma]^{(2)}_x$ the elements of $D^{(2)}$ such that $\tau_D(\gamma(0))=x.$

Consider the  dual bundle $\pi_D:D^*\rightarrow Q$. If $\beta\in \Gamma(\pi_D),$ then
we define the \emph{linear function} $\hat{\beta}: D\rightarrow {\mathbb R}$ by
\[
\hat{\beta}(v_x)=\langle \beta(x), v_x\rangle,
\]
for all $v_x\in D$.

We define the complete lift of a section in an analogous way to proposition (\ref{prop:levantcpleto}) as following.

\begin{definition}\label{propo1-A}
The \emph{complete lift} $X^c$ of a section $X\in\Gamma(\tau_D)$ is the unique vector field $X^c\in\mathfrak{X}(D)$ which verifies the following two conditions:
\begin{enumerate}
\item $X^c$ is projectable over $X$ by means of $(\tau_{D})_*$; that is, $(\tau_{D})_*X^c=X$.
\item $X^c(\hat{\alpha})=\widehat{{\mathcal L}_X\alpha}$, for all $\alpha\in \Gamma(\pi_D)$.
\end{enumerate}
\end{definition}
Here, ${\mathcal L}_X\alpha\in \Gamma (\tau_D)$ denotes the Lie derivative of $\alpha\in\Gamma(\pi_D)$   with respect to $X\in\Gamma(\tau_D)$
(see \cite{dLMMalgebroides} for details):
\[
\langle {\mathcal L}_X\alpha, Y\rangle=\rho_D(X)(\langle \alpha, Y\rangle) - \langle \alpha,  \lcf X,Y\rcf_D\rangle,\,\  \forall\ Y\in \Gamma(\tau_D).
\]

 Let us assume that we
have fixed coordinates $(q^i)$ in $Q$ and an arbitrary frame $\{e_A\},$ then an arbitrary section $X\in \Gamma(\tau_D)$ will have an expression
$X=\tilde{X}^A(q) e_A.$
Moreover, the new frame induces a new system of coordinates $(q^i, y^A)$ on $D$, where $v_x=y^A e_A(x)$ for any $v_x\in D$.
Using Proposition \ref{propo1-A} or by a direct change of coordinates, it is not hard to prove that the complete lift $X^c$ can be rewritten as
\[
X^c=\rho_A^i\tilde{X}^A\frac{\partial}{\partial q^i}+\left(\rho^i_B\frac{\partial \tilde{X}^C}{\partial q^i}- \mathcal{C}_{AB}^C
\tilde{X}^A \right) y^B\frac{\partial}{\partial y^C}\; .
\]

Another notion to be used later is that of vertical lift.
\begin{definition}
   The \emph{vertical lift} $X^v$ of a section $X$ of $D$ is the unique vector field $X^v\in\mathfrak{X}(D)$ verifying the following conditions:
   \begin{enumerate}
      \item $X^v(\tau^*_D f)=0,$ for all $f\in C^{\infty}(Q),$
      \item $X^v(\hat\alpha)=\tau^*_D(\langle\alpha,X\rangle),$ for all $\alpha\in\Gamma(\pi_D).$
   \end{enumerate}
\end{definition}

In coordinates, $\displaystyle X^v=\tilde{X}^A\frac{\partial}{\partial y^A}.$

For $I\subseteq {\mathbb R},$ a \emph{time-dependent section} $X$ is a smooth mapping
$
X: I\times Q\rightarrow D$ such that $X(t, x)\in D_x$.

\begin{definition}
The \emph{tangent lift} $X^T$ of a time-dependent section on $Q$ is the unique time-dependent vector field $X^T\in\mathfrak{X}(pr_D)$, where $pr_D:I\times D\rightarrow D$, verifying the following two conditions:
\begin{enumerate}
\item $X^T$ is projectable over $\rho_D(X)$ by means of $\tau_D$; that is, $(\tau_D)_*X^T=\rho_D(X).$
\item $(X^T)_{t,v_x}(\hat{\alpha})=\widehat{{\mathcal L}_{X_t}\alpha}(v_x)+\frac{d}{dt}\langle \alpha_x, X_x(t)\rangle$, for all $\alpha\in \Gamma(\pi_D)$.
\end{enumerate}
\end{definition}
Here, $X(t,x)=X_t(x)=X_x(t).$

We have the diagram

\begin{equation*}
\xymatrix{
I\times D\ar[d]_{(id_{\mathbb{R}},\tau_D)}\ar[rr]^{X^T}&&TD\ar[d]^{\tau_{TD}}\\
I\times Q\ar[rr]^{X}&&D
}
\end{equation*}
If  $X=\tilde{X}^A(t,q)e_A,$ then  we have
\[
X^T=\rho^{i}_A(q)\tilde{X}^A(t,q)\frac{\partial}{\partial q^i}+\left[\frac{\partial\tilde{X}^C}
{\partial t}(t,q)+\left( \rho^{i}_B(q)\frac{\partial \tilde{X}^C}{\partial q^i}(t,q)- \mathcal{C}_{AB}^C(q)\tilde{X}^A(q)\right) y^B\right]\frac{\partial}{\partial y^C}.
\]
Similarly, we can introduce the \emph{vertical lift} $X^V\in\mathfrak{X}(pr_D)$ of a time-dependent vector field $X$ as
\[
X^V(t, v_x)=(X_t)^v(v_x),
\]
where $X_t$ is the
vector field on $Q$ defined by $X_t(x) = X(t, x)$.

\vspace{.2cm}
The following definition will be useful in the sequel.
\begin{definition}\label{def:E-L-Oper-A}
The \emph{Euler-Lagrange operator} associated with a 1-form $\mu\in \Gamma (\pi_{TD})\equiv \Lambda^1D$  is the mapping ${\mathcal E}_{\mu}: D^{(2)}\rightarrow D^*$ defined by
\[
 \left\langle {\mathcal E}_{\mu} ([\gamma]^{(2)}_x), X(x)\right\rangle= \frac{d}{dt}\langle \mu,  X^v\rangle ([\gamma]^{(2)}_x)-\langle \mu, X^c\rangle (\gamma(x)),
\]
 for any $X\in {\Gamma}(\tau_D)$.
\end{definition}

Observe that, if $X$ is time-dependent, we have that
\begin{equation}\label{pok-A}
 \left\langle {\mathcal E}_{\mu} ([\gamma]^{(2)}_x), X(t,x)\right\rangle= \frac{d}{dt}\langle \mu,  X^V\rangle(t, [\gamma]^{(2)}_x)-\langle \mu, X^T\rangle(t, \gamma(x)).
 \end{equation}

For a function $L: D\rightarrow \mathbb{R}$  we have
\begin{equation}\label{EulerLagnonholframe-A}
 \langle {\mathcal E}_{\mathrm{d}L} , X\rangle= \tilde{X}^A\left[
 \frac{d}{dt}\left(\frac{\partial L}{\partial y^A}\right)
 -\rho^i_A\frac{\partial L}{\partial q^i}+{\mathcal C}_{AB}^C
y^B\frac{\partial L}{\partial y^C}\right].
\end{equation}

\subsection{Generalized variational problem on skew-symmetric algebroids.}
As in Section \ref{sect:GVCtangbundle}, we can directly define generalized variational calculus on skew-symmetric algebroids.
\begin{definition}
A {\bf generalized variational problem} on a skew-symmetric algebroid $D$ is determined by a submanifold $\Sigma$ of $T^*D.$
\end{definition}

\begin{definition}
A  {\bf solution} of the generalized variational problem determined by $\Sigma\subset T^*D$ is an admissible curve
$\gamma:I\subseteq\mathbb{R}\rightarrow D$ such that
there exists a curve $\mu:I\rightarrow\Sigma$  verifying $\pi_{TD}(\mu(t))=\gamma(t)$ and for all time-dependent section $X$ of $\tau_D:D\rightarrow Q,$
\begin{equation}\label{pou-A}
\int_I\left\langle\mu_{\gamma(t)},X^T(t,\gamma(t))\right\rangle \,\mathrm{d}t=0.
\end{equation}
%\end{enumerate}
\end{definition}

Taking the canonical projection $\pi_{TD}: T^*D\rightarrow D,$ we define the subset $C= \pi_{TD}(\Sigma)$ (\emph{kinematical constraints}) and we have
the following diagram:

\begin{equation*}
\xymatrix{
&\Sigma\ar[d]_{(\pi_{D})_{|\Sigma}}\ar@{^{(}->}[rr]^{i_{\Sigma}}&&T^*D\ar[d]^{\pi_{D}}\\
&C\ar@{^{(}->}[rr]&&D\ar[d]_{\tau_D}\\
I\ar[rrr]_{\sigma=\tau_D\circ\gamma}\ar[ur]^{\gamma}\ar@/^2.5pc/[uur]^{\mu} & &&Q
}
\end{equation*}

\vspace{.2cm}
Since $\pi_{D}(\mu(t))=\gamma(t)$ then $\text{Im}(\gamma)\subseteq C.$ Moreover, if $\sigma=\tau_D\circ\gamma,$ since $\gamma$ is admissible then
$\dot{\sigma}(t)=\rho_D(\gamma(t)).$

Analogously to Section \ref{sect:GVCtangbundle}, we deduce that an admissible curve $\gamma:I\rightarrow D$ is a solution of the generalized variational problem if there exists $\mu:I\rightarrow\Sigma\subseteq T^*D$ verifying $$\mathcal{E}_{\mu}([\gamma]^{(2)}_x)=0.$$
 In  local coordinates, we assume that $\Sigma$ is determined by the vanishing of constraints $\Phi^{\alpha}=0$ where $\Phi^{\alpha}: T^*D\rightarrow \mathbb{R}$.
 A curve $\gamma:I\rightarrow D,$ locally given by  $\gamma(t)=(q^i(t),y^A(t)),$ is admissible if $\rho_A^i(q(t))y^A(t)=\mathrm{d}q^i(t)/\mathrm{d}t$.
Therefore, we seek a
curve $\mu:I\rightarrow T^*D$ over $\gamma,$ locally expressed as  $\mu=(\mu_i(\gamma(t))\mathrm{d}q^i+\tilde{\mu}_A(\gamma(t))\mathrm{d}y^A)\mid_{\gamma},$
such that $\mathcal{E}_{\mu}([\gamma]^{(2)}_x)=0;$ that is
$$\frac{d}{dt}\tilde{\mu}_A-\rho_A^j\mu_j+\mathcal{C}^C_{AB}y^B\tilde{\mu}_C=0,$$
or equivalently, $$\frac{\partial\tilde{\mu}_A}{\partial q^i}\rho_B^iy^B+\frac{\partial\tilde{\mu}_A}{\partial y^B}\dot{y}^B-\rho^j_A\mu_j+\mathcal{C}^C_{AB}y^B\tilde{\mu}_C=0\;.
$$
Summarizing, we have the following set of equations:
\begin{equation}\label{ec:CVG-A}
\left\{
\begin{array}{rl}
\displaystyle\frac{d}{dt}\tilde{\mu}_A-\rho^j_A\mu_j+\mathcal{C}^C_{AB}y^B\tilde{\mu}_C&=0 \\
\\
\displaystyle \Phi^{\alpha}(q^i,y^A,\mu_i,\tilde{\mu}_A)&=0
\\
\displaystyle\frac{{d}q^i}{{d}t}&=\rho^i_A y^A.
\end{array}
\right.
\end{equation}
In the sequel, we will describe some particular examples of generalized variational calculus on skew-symmetric algebroids.

\subsection{Lagrangian mechanics on skew-symmetric algebroids}
 Given a function \newline $L:D\rightarrow\mathbb{R}$, we take
$\Sigma=Im(\mathrm{d}L)=\mathrm{d}L(D)\subseteq T^*D$. In this case, $C=D$ and we try to find an admissible curve $\xi:I\rightarrow D$ such that
$\displaystyle\int_{t_0}^{t_1}\left\langle \mathrm{d}L(\xi(t)),X^T(t, \xi(t))\right\rangle\, \mathrm{d}t=0,$ for all time dependent section $X$ of
$\tau_D:D\rightarrow Q.$ From this equation we derive the Euler-Lagrange equations (see \cite{martinez01,weinstein78}) given by
\begin{equation}\label{EcLagSkew}
\left\{
\begin{array}{rl}
\displaystyle\frac{d}{dt}\left(\frac{\partial L}{\partial y^A}\right)
 -\rho^i_A\frac{\partial L}{\partial q^i}+{\mathcal C}_{AB}^C
y^B\frac{\partial L}{\partial y^C}&=0 \\
\\
\dot{q}^i&=\rho^i_A y^A.
\end{array}
\right.
\end{equation}

\subsubsection{The Euler-Poincar\'e equations} See \cite{Holm1, Holm2, marsden3}).
In this case, we have a Lagrangian $l:\mathfrak{g}\rightarrow\mathbb{R}$ defined on the Lie algebra $\mathfrak{g}$ of a Lie group $G$ and we consider $\Sigma=\mathrm{d}l(\mathfrak{g})\subseteq T^*\mathfrak{g}\simeq\mathfrak{g}\times\mathfrak{g}^*.$
A time-dependent section is then a curve $\eta:I\rightarrow\mathfrak{g}$  and, therefore, its tangent lift is the time-dependent vector field on ${\mathfrak g}$
defined by  $\eta^T(t,\xi)=(\xi,\dot{\eta}(t)+\text{ad}_{\xi}\eta(t))\in\mathfrak{g}\times\mathfrak{g}\simeq T\mathfrak{g}.$

Hence, a curve $\xi:I\rightarrow\mathfrak{g}$ is a solution of the generalized variational problem if

\begin{eqnarray*}
   \displaystyle 0&=&\int_{t_0}^{t_1}\left\langle\mathrm{d}l(\xi(t)),\eta^T(t,\xi(t))\right\rangle\mathrm{d}t\\
   &=& \displaystyle\int_{t_0}^{t_1}\left\langle\mathrm{d}l(\xi(t)),\dot{\eta}(t)+\text{ad}_{\xi(t)}\eta(t)\right\rangle
   \mathrm{d}t\\
   &=&   \displaystyle\int_{t_0}^{t_1}\left\langle
   \frac{\partial l}{\partial\xi}(\xi(t)),
   \dot{\eta}(t)\right\rangle+\left\langle\text{ad}^*_{\xi(t)}
   \frac{\partial l}{\partial\xi}(\gamma(t)),\eta(t)\right\rangle\mathrm{d}t\\
   &=&   \displaystyle\int_{t_0}^{t_1}\left\langle\left(\text{ad}^*_{\xi}
   \frac{\partial l}{\partial\xi}-\frac{\mathrm{d}}{\mathrm{d}t}
   \frac{\partial l}{\partial\xi}\right)\Big |_{\xi(t)},\eta(t)\right\rangle\mathrm{d}t.
\end{eqnarray*}

From this,  we deduce the classical Euler-Poincar\'e equations $$\frac{{d}}{{d}t}\left(\frac{\partial l}{\partial\xi}\right)=\text{ad}^*_{\xi}\frac{\partial l}{\partial\xi}.$$

\subsubsection{Nonholonomic mechanics.}
Let $\mathcal{L}:TQ\rightarrow\mathbb{R}$ be a mechanical Lagrangian on $TQ,$ that is, $\displaystyle \mathcal{L}(v_q)=1/2\;g(v_q,v_q)-V(q)$ where $g$ is a Riemannian metric on $Q$ and $V:Q\rightarrow\mathbb{R}$ is a potential function.

Additionally, in the case of nonholonomic mechanics, we have a regular distribution $D\subseteq TQ.$ Using the Riemannian metric $g,$ we consider the
Riemannian orthogonal decomposition $TQ=D\oplus D^{\perp, g}$ and we denote by $i_D:D\hookrightarrow TQ$ the canonical inclusion and by
$P:TQ\rightarrow D$ the associated orthogonal projector. We induce a skew-symmetric algebroid structure $\displaystyle \lcf X,Y\rcf_D=P[i_D X,i_D Y],$
for $X,Y\in\Gamma(\tau_D)$ (See for instance \cite{MaríaHamJac,MR2492630}). Observe that the bracket $\lcf .,. \rcf_D$ does not satisfy the Jacobi identity in general, and therefore $(D, \lcf .,.\rcf_D,\rho_D)$ is a skew-symmetric algebroid where $\rho_D=i_D: D\hookrightarrow TQ$.

\vspace{.2cm}
Let $L:D\rightarrow\mathbb{R}$ be the restriction of the Lagrangian $\mathcal{L}$ to $D$,  i.e., $L=\mathcal{L}\circ i_D.$

We are now able to apply our generalized variational calculus to the mechanical system determined by $L:D\rightarrow\mathbb{R}$ and $D,$
with its mentioned skew-symmetric algebroid structure induced from the orthogonal projection of the standard Lie bracket to $D$.
Hence, a solution of the nonholonomic problem is a curve $\gamma:I\rightarrow D$ such that $\gamma$ is admissible and there exists a curve $\mu:I\rightarrow Im(\mathrm{d}L)\subseteq T^*D$ such that $\pi_{TD}(\mu(t))=\gamma(t)$ and, for all time-dependent section $X$ of $\tau_D:D\rightarrow Q,$
$$\int_I\left\langle\mu_{\gamma(t)},X^T(t,\gamma(t))\right\rangle\mathrm{d}t=0,$$
where $X^T$ is the tangent lift given by the induced skew-symmetric algebroid structure. Then, the equations of the nonholonomic problem are equations (\ref{EcLagSkew}) which are the Lagrange-d'Alembert's equations in this context (see \cite{MaríaHamJac}).
It is easy to adapt the previous calculations to nonholonomic systems with symmetries (see \cite{MR2542960,MR2492630}).

\subsection{Hamiltonian mechanics on skew-symmetric Lie algebroids}\label{mecHam-A}

Let $H:D^*\rightarrow\mathbb{R}$ be a function where $\pi_D:D^*\rightarrow Q$ is the dual bundle of an skew-symmetric algebroid $\tau_D:D\rightarrow Q.$
In a similar way to section \ref{mecHam}, it is defined an antisymplectomorphism $\mathcal{R}:T^*D^*\rightarrow T^*D$ (see \cite{GraGra08,MackPX94}).
In local coordinates, if $(q^i)$ are coordinates on $Q$ and $\{e_A\}$ is a basis of sections of $\tau_D:D\rightarrow Q,$ then we have the dual basis of section $\{e^A\}$ on $D^*$ (that is, $e^A(e_B)=\delta_B^A$). This dual basis induces coordinates $(q^i,p_A)$ on $D^*.$
The antisymplectomorphism $\mathcal{R}$ is given by $\mathcal{R}(q^i,p_A,\gamma_i,\gamma^A)=(q^i,\gamma^A,-\gamma_i,p_A).$

Now, we construct the submanifold $\Sigma_H\subseteq T^*D$ by $\Sigma=\mathcal{R}(\mathrm{d}H(D^*)).$
Locally, $$\displaystyle\Sigma_H=\left\{(q^i,y^A,\mu_i,\tilde{\mu}_A)\; \mid \; y^A=\frac{\partial H}{\partial p_A}(q,\tilde{\mu}),\mu_i=-\frac{\partial H}{\partial q^i}(q,\tilde{\mu})\right\}.$$

An admissible curve $\gamma:I\rightarrow D$ is a solution of the Hamiltonian problem given by $H:D^*\rightarrow\mathbb{R}$ if there exists a curve $\mu:I\rightarrow\Sigma$ such that
$$\int_I\langle\mu(t),X^T(t,\gamma(t))\rangle\mathrm{d}t=0,$$ for all time-dependent section $X$ of $\tau_D:D\rightarrow Q.$

If $\mu:I\rightarrow\Sigma$ is given by $\mu(t)=(q^i(t),y^A(t),\mu_i(t),\tilde{\mu}_A(t)),$ then equations (\ref{ec:CVG-A}) are equivalent to the following set of equations:
\begin{eqnarray*}
   \displaystyle\frac{{d\tilde{\mu}_A}}{{d}t}+\rho_A^j\frac{\partial H}{\partial q^j}+\mathcal{C}_{AB}^Cy^B\tilde{\mu}_C&=&0\\
   \displaystyle y^A&=&\frac{\partial H}{\partial p_A}(q,\tilde{\mu})\\
   \dot{q}^i&=&\rho_A^iy^A
\end{eqnarray*}

or, equivalently,

\begin{eqnarray}
   \displaystyle\frac{{d\tilde{\mu}_A}}{{d}t}&=&-\rho_A^j\frac{\partial H}{\partial q^j}-\mathcal{C}_{AB}^C\frac{\partial H}{\partial p_B}
\tilde{\mu}_C\\
   \displaystyle \frac{{d}q^i}{{d}t}&=&\rho^i_A\frac{\partial H}{\partial p_A},
\end{eqnarray}
which are Hamilton's equations in the context of skew-symmetric Lie algebroids (see \cite{dLMMalgebroides}).

\subsubsection{Lie Poisson equations}
Let ${\mathfrak g}$ be a Lie algebra, ${\mathfrak g}^*$ its dual and
$H:\mathfrak{g}^*\rightarrow\mathbb{R}$ a hamiltonian function.  We construct $\displaystyle\Sigma_H=\mathcal{R}(\mathrm{d}H(\mathfrak{g}^*))=\left\{(\xi,\alpha)\in\mathfrak{g}\times\mathfrak{g}^*
\equiv T^*\mathfrak{g}\; |\; \xi=\frac{\partial H}{\partial\alpha}\right\}.$

A solution of the generalized variational problem is characterized by
$$0 = \int_{t_0}^{t_1}\left\langle\alpha,\dot{\eta}+\text{ad}_{\frac{\partial H}{\partial \alpha}}\eta\right\rangle\mathrm{d}t=
\int_{t_0}^{t_1}\left\langle\text{ad}^*_{\frac{\partial H}{\partial\alpha}}\alpha-\frac{\mathrm{d}\alpha}{\mathrm{d}t},\eta\right\rangle\mathrm{d}t$$
which give us the classical \emph{Lie-Poisson equations} (see \cite{Holm1,Holm2}) $$\frac{\mathrm{d}\alpha}{\mathrm{d}t}=\text{ad}^*_{\frac{\partial H}{\partial\alpha}}\alpha.$$

\subsection{Constrained variational calculus}\label{calcvariacrestringido-A}

Now, we study the case of variational cons\-trained calculus on the setting of skew-symmetric algebroids (see \cite{GraGra08,DianaJCDavides07}). We will see how to define a submanifold of $T^*D$ to apply our generalized variational calculus and to derive the corresponding equations in this case.

The variational constrained problems are determined by a pair $(C, l)$ where $C$ is a submanifold of $D$, with inclusion $i_C: C\hookrightarrow D$, and $l: C\rightarrow \mathbb{R}$ is a Lagrangian function defined \emph{only} along $C$.
We will consider
 \begin{multline*}
    \Sigma _l = \bigl\{ \mu \in T ^\ast D \mid \pi _{TD} (\mu) \in C \text{
        and } \left\langle \mu, v \right\rangle = \left\langle
        \mathrm{d} l , v \right\rangle, \\
      \text{ for all } v \in T C \subseteq T D \text{ such that } \tau
      _{TD} (v) = \pi _{TD} (\mu) \bigr\}.
  \end{multline*}

We can also write $\Sigma_l$ as
\begin{equation*}
\Sigma_l=\{\mu \in T ^\ast D\mid i_C^*\mu=dl \}=(\mathrm{d}L+\nu^{*}(C))\mid_{C}.
\end{equation*}
Here $L: D\rightarrow \mathbb{R}$ is an arbitrary extension of $l$ to $D$ (i.e., $l\circ i_C=L$) and
$\nu^{*}(C)$  is the conormal bundle of $C$. Considering a curve $\mu: [t_0,t_1]\rightarrow \Sigma_l$ as $\mu(t)=\mathrm{d}L(\gamma(t))+\nu(t)$ where
$\nu(t)\in \left(\nu^{*}(C)\right)\mid_{\gamma(t)},$ and $\gamma$ is an admissible curve satisfying $\gamma(t)\in C\subseteq D,$
then
\begin{equation*}
 \displaystyle \int_{t_0}^{t_1}\left\langle\mu(t),X^T(t, \dot{\sigma}(t))\right\rangle\ \mathrm{d}t=0
\end{equation*}
becomes
\begin{eqnarray*}\label{ecvak}
{\mathcal E}_{{\mathrm{d}L+\nu}}&=&0 \label{plm1-A}\\
\gamma(t)&\in& C,\quad \forall t\in[t_0,t_1], \label{plm2-A}
\end{eqnarray*}
and a solution is a pair $(\gamma, \nu)$ with $\gamma: I\rightarrow D$ admissible  and $\nu(t)\in \left[\nu^{*}(C)\right]_{\gamma(t)}$.

If we assume that we have fixed local constraints which determine $C$ by their vanishing, i.e.,  $\phi^{\alpha}(q, y)=0$, $1\leq\alpha\leq n$ where
$n=\hbox{dim } TD-\hbox{dim } C$, then
\[
\left[\nu^{*}(C)\right]\mid_{\dot{\sigma}(t)}=\hbox{span }\{ \mathrm{d}\phi^{\alpha}(\dot{\sigma}(t))\}
\]
and $\nu(t)=\lambda_{\alpha}(t)\mathrm{d}\phi^{\alpha}(\dot{\sigma(}t))$ for some Lagrange multiplier $\lambda_{\alpha}.$ So, equations (\ref{ecvak}) are
\begin{equation*}
\left\{
 \begin{array}{rcl}
{\mathcal E}_{\mathrm{d}L+\lambda_{\alpha} d\phi^{\alpha}} &=&0\\
\phi^{\alpha}(\gamma(t))&=&0.
\end{array}
\right.
\end{equation*}

As in Section \ref{calcvariacrestringido}, we derive the following system of equations
\begin{equation*}
\left\{
 \begin{array}{rl}
\dot{q}^i-\rho^i_Ay^A &=0\\
\\
\displaystyle \frac{d}{dt}\left(\frac{\partial(L+\lambda_{\alpha}\Phi^{\alpha})}{\partial y^A}\right)-\rho^i_A\frac{\partial(L+\lambda_{\alpha}\phi^{\alpha})}{\partial q^i}+\mathcal{C}^C_{AB}y^B\frac{\partial(L+\lambda_{\alpha}\phi^{\alpha})}{\partial y^C}&=0\\
\\
\phi^{\alpha}(q^i,y^A)&=0
\end{array}
\right.
\end{equation*}
which is a generalization of equations (\ref{gralvakonnoncanon}).

%If we describe the problem defined by $\Sigma_l=\{\mu\in T^*D:i^*_C\mu=\mathrm{d}l\}$ where $C$ is now locally described  by the constraints   $y^{\alpha}=\phi^{\alpha}(q^i,y^a)%$ with $1\leq a\leq m,$ $m+1\leq\alpha\leq rank D$ and $\dim C=\dim Q+m.$ Again proceeding as in section \ref{calcvariacrestringido} we obtain exactly the equations %
%(\ref{ecuacvakonnoncanon}) where now $1\leq a\leq m,$ $m+1\leq\alpha\leq rank D$ and $1\leq i\leq\dim Q.$

\subsection{Optimal Control Theory}
An \emph{optimal control problem} on a skew-symmetric algebroid is given by a quadruple  $({C}, Q,\Gamma, { L})$ where
$\tau_{{ C},Q}\colon {C} \rightarrow Q$ is the control bundle, $\Gamma$ is a vector field defined along $\tau_{{ C},Q}$ and   ${ L} \colon { C} \rightarrow \mathbb{R}$ is a cost function whose associated functional
must be minimized.
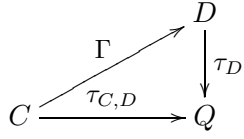
\begin{equation*}\xymatrix{&& D \ar[d]^{\txt{\small{$\tau_D$}}} \\ {C}
\ar[rr]^{\txt{\small{$\tau_{{C},D}$}}} \ar[rru]^{\txt{\small{$\Gamma$}}}& & Q}
\end{equation*}
Locally, we have that $y^A=\Gamma^A(q,u)$.
From the optimal control data  $( C, Q, \Gamma, { L})$ we construct Pontryagin's hamiltonian $H:D^*\times_Q  {\mathcal C}\longrightarrow \mathbb{R}$ given by
\begin{equation}\label{pont-A}
H(\alpha_q, u_q)=\langle \alpha_q, \Gamma(u_q)\rangle-L(u_q),
\end{equation}
where $u_q\in C_q$ and $\alpha_q\in D_q^*$.
In coordinates, $H(q^A, p_A, u^a)=p_A\Gamma^A(q, u)-L(q, u)$.

Also the optimal control solutions can be characterized using generalized variational calculus. We define the subset
\[
\Sigma=\{\mu\in T^*D\; \mid\; \Gamma^*\mu=\mathrm{d}L\}\, .
\]
If we take an arbitrary element $\mu=\mu_i\mathrm{d}q^i+\tilde\mu_A\mathrm{d}y^A$ in $T^*D,$ then a solution  curve for the generalized variation calculus associated to $\Sigma$ is given by the following system of equations
\begin{equation}\label{H-PPrinciple-A}
   \left\{
   \begin{array}{rcl}
      \displaystyle \frac{d}{dt}\tilde{\mu}_A=\left(\frac{\partial L}{\partial q^i}-\tilde{\mu}_A\frac{\partial \Gamma^A}{\partial q^i}\right)\rho_A^i-\mathcal{C}_{AB}^Cy^B\tilde{\mu}_C\\
      \displaystyle \tilde\mu_A\frac{\partial \Gamma^A}{\partial u^{a}}=\frac{\partial L}{\partial u^a}\\
      \\
      \dot{q}^i=\rho_A^i\Gamma^A
   \end{array}
   \right.
\end{equation}
for some curve $u(t)=(u^a(t)).$
Alternatively, in terms of Pontryagin's Hamiltonian $H$ the equations are rewritten as follows:
\begin{equation}
   \left\{
   \begin{array}{rcl}
      \displaystyle \frac{d}{dt}\tilde{\mu}_A=-\frac{\partial H}{\partial q^i}\rho_A^i-\mathcal{C}_{AB}^C\Gamma^B\tilde{\mu}_C\\
      \displaystyle \frac{\partial H}{\partial u^{a}}=0\\
      \displaystyle\dot{q}^i=\rho_A^i\Gamma^A.
   \end{array}
   \right.
\end{equation}
See these equations in references \cite{DianaJCDavides07,MR2050148}.

\section{\textbf{\textsc{Discrete generalized variational calculus on $Q\times Q$}}}\label{DGVCenQxQ}

In this section, we will develop a discrete version of the generalized variational calculus. For that, we will only need to have a subset of an appropriate cotangent bundle  and to introduce the notions of discrete complete and vertical lifts.

The main motivation will be the derivation of numerical integrators for the co\-rres\-pon\-ding continuous systems which preserve some of their geometric or invariance properties, see \cite{marsden-west}.

\subsection{Discrete geometry}

The discrete notion of the tangent bundle $TQ$ is the cartesian product of two copies of $Q,$ that is, $Q\times Q.$ Now, we have two canonical projections $\alpha:Q\times Q\rightarrow Q$ defined by $\alpha(q,\tilde{q})= q$ and $\beta:Q\times Q\rightarrow Q$ defined by  $\beta(q,\tilde{q})=\tilde{q},$ where $q, \tilde{q}\in Q$.

Given a curve $\sigma:\mathbb{Z}\rightarrow Q,$ we define its tangent lift $\tilde{\sigma}:\mathbb{Z}\rightarrow Q\times Q$ as follows: $$\tilde{\sigma}(k)=(\sigma(k),\sigma(k+1)),$$ for all $k\in\mathbb{Z}.$

As in the continuous case, we need to introduce the notion of \emph{discrete complete lift} $X^c\in\mathfrak{X}(Q\times Q)$ of a vector field $X\in\mathfrak{X}(Q).$ It is defined by $$X^c(q,\tilde{q})=(X(q),X(\tilde{q}))\in T_qQ\times T_{\tilde{q}}Q\equiv T_{(q, \tilde{q})} (Q\times Q)\; .$$
Moreover, we have two notions of \emph{discrete vertical lifts} of $X$ given by the following formulas:
$$X^{v_{\alpha}}(q,\tilde{q})=(X(q),0_{\tilde{q}})\quad\text{and}\quad X^{v_{\beta}}(q,\tilde{q})=(0_q,X(\tilde{q})).$$

In the same way as in section \ref{sect:fibtangente}, for all $\mathbb{Z}$-dependent vector field $X:\mathbb{Z}\times Q\rightarrow TQ,$ we have
 its \emph{discrete tangent lift} $X^T:\mathbb{Z}\times Q\times Q\rightarrow TQ\times TQ$ defined by $$X^T(k,q,\tilde{q})=(X(k,q),X(k+1,\tilde{q})),$$
and we have $X^{V_{\alpha}}(k,q,\tilde{q})=(X_k)^{v_{\alpha}}(q,\tilde{q})$ and $X^{V_{\beta}}(k,q,\tilde{q})=(X_k)^{v_{\beta}}(q,\tilde{q}),$ where $X_k(q,\tilde{q})=X(k,q,\tilde{q}).$

\vspace{.3cm}
\begin{definition}
   The  \emph{discrete Euler-Lagrange operator} associated to a 1-form $\mu\in\Gamma(\pi_{T(Q\times Q)})\equiv\Lambda^1(Q\times Q) $ is the mapping $\mathcal{E}_{\mu}^d:Q\times Q\times Q\rightarrow T^*Q$ defined by
   $$\left\langle\mathcal{E}_{\mu}^d(q,\tilde{q},\overline{q}),X(\tilde{q})\right\rangle=
   \langle\mu_{(\tilde{q},\overline{q})},X^{v_{\alpha}}(\tilde{q},\overline{q})\rangle-
   \langle\mu_{(q,\tilde{q})},X^{v_{\alpha}}(q,\tilde{q})\rangle+
   \langle\mu_{(q,\tilde{q})},X^c(q,\tilde{q})\rangle$$
   or, alternatively,
   $$\left\langle\mathcal{E}_{\mu}^d(q,\tilde{q},\overline{q}),X(\tilde{q})\right\rangle=-[
   \langle\mu_{(\tilde{q},\overline{q})},X^{v_{\beta}}(\tilde{q},\overline{q})\rangle-
   \langle\mu_{(q,\tilde{q})},X^{v_{\beta}}(q,\tilde{q})\rangle-
   \langle\mu_{(\tilde{q},\overline{q})},X^c(\tilde{q},\overline{q})\rangle],$$
   for all $X\in\mathfrak{X}(Q).$
\end{definition}

\vspace{.2cm}
If $\mu\in\Lambda^1(Q\times Q),$ then we can decompose $\mu(q,\tilde{q})=(\mu_1(q,\tilde{q}),\mu_2(q,\tilde{q}))$ where $\mu_1(q,\tilde{q})\in T_q^*Q$ and $\mu_2(q,\tilde{q})\in T^*_{\tilde{q}}Q,$ then $$\mathcal{E}^d_{\mu}(q,\tilde{q},\overline{q})=\mu_1(\tilde{q},\overline{q})+\mu_2(q,\tilde{q})\in T^*_{\tilde{q}}Q.$$

In particular, for a function $L_d:Q\times Q\rightarrow\mathbb{R}$ we have $$\mathcal{E}_{\mathrm{d}L_d}^d(q,\tilde{q},\overline{q})=D_1L_d(\tilde{q},\overline{q})+D_2L_d(q,\tilde{q})$$
where we use the notation $D_1$ and $D_2$ for the decomposition
$$\displaystyle \mathrm{d}L_d(q,\tilde{q})=(D_1L_d(q,\tilde{q}),D_2L_d(q,\tilde{q})),$$ with $D_1L_d(q,\tilde{q})\in T^*_qQ$ and $D_2L_d(q,\tilde{q})\in T^*_{\tilde{q}}Q.$

Now, if $X:\mathbb{Z}\times Q\rightarrow TQ$ with $X(k,q)\in T_qQ,$ for all $k\in\mathbb{Z}$ and for all $q\in Q,$ then

\begin{eqnarray}\label{E-Ldiscreto}
\langle\mathcal{E}_{\mu}^d(q,\tilde{q},\overline{q}),X(k,\tilde{q})\rangle &=&
\langle\mu_{(\tilde{q},\overline{q})},X^{V_{\alpha}}(k,\tilde{q},\overline{q})\rangle-
\langle\mu_{(q,\tilde{q})},X^{V_{\alpha}}(k-1,q,\tilde{q})\rangle \nonumber \\
& &+\langle\mu_{(q,\tilde{q})},X^T(k-1,q,\tilde{q})\rangle \nonumber \\
&= &-(\langle\mu_{(\tilde{q},\overline{q})},X^{V_{\beta}}(k+1,\tilde{q},\overline{q})\rangle-
\langle\mu_{(q,\tilde{q})},X^{V_{\beta}}(k,q,\tilde{q})\rangle \\
& &-\langle\mu_{(\tilde{q},\overline{q})},X^T(k,\tilde{q},\overline{q})\rangle). \nonumber
\end{eqnarray}

\subsection{Discrete generalized variational problem}

With the above definitions, we can introduce the notion of generalized variational calculus in the context of discrete mechanics as follows.

\begin{definition}
A \textbf{discrete generalized variational problem} is determined by a submanifold $\Sigma_d\subseteq T^*(Q\times Q)\equiv T^*Q\times T^*Q.$
\end{definition}

\begin{definition}\label{discretemechanicsSolution}
A \textbf{solution} of the discrete generalized variational problem determined by $\Sigma_d\subseteq T^*(Q\times Q)$ is a curve $\sigma:\mathbb{Z}\rightarrow Q$ such that there exists a curve $\mu:\mathbb{Z}\rightarrow\Sigma_d$ such that $\pi_{T(Q\times Q)}(\mu(k))=\pi_{T(Q\times Q)}(\mu_{(\sigma(k),\sigma(k+1))})=\tilde{\sigma}(k),$ where $\tilde{\sigma}(k)=(\sigma(k),\sigma(k+1))\in Q\times Q$ and, for all $\mathbb{Z}$-dependent section $X:\mathbb{Z}\times Q\rightarrow TQ,$
$$\sum_{k=0}^{N-1}\left\langle\mu_{(\sigma(k),\sigma(k+1))},X^T(k,\sigma(k),\sigma(k+1))\right\rangle=0.$$ \end{definition}

Denoting by $C_d=\pi_{T(Q\times Q)}(\Sigma_d)$ the set of \emph{discrete kinematical constraints}, we have the diagram

\begin{equation*}
\xymatrix{
&\Sigma_d\ar[d]_{(\pi_{T(Q\times Q)})_{|\Sigma_d}}\ar@{^{(}->}[rr]^{i_{\Sigma_d}}&&T^*(Q\times Q)\ar[d]^{\pi_{T(Q\times Q)}}\\
&C_d\ar@{^{(}->}[rr]&&Q\times Q\ar@<1ex>[d]_{\alpha\ }\ar[d]^{\ \beta}\\
\mathbb{Z}\ar[rrr]_{\sigma}\ar[ur]^{\tilde{\sigma}}\ar@/^4.5pc/[uur]^{\mu} & &&Q
}
\end{equation*}

Since $\displaystyle\mu_{\tilde{\sigma}(k)}\in T^*_{\tilde{\sigma}(k)}(Q\times Q)\simeq T^*_{(\sigma(k),\sigma(k+1))}(Q\times Q)=T^*_{\sigma(k)}Q\times T^*_{\sigma(k+1)}Q,$ then we can write $\displaystyle\mu_{\tilde{\sigma}(k)}=((\mu_1)_{{\sigma}(k)},(\mu_2)_{{\sigma}(k+1)})\in T^*_{\sigma(k)}Q\times T^*_{\sigma(k+1)}Q.$

\vspace{.2cm}
Observe that from (\ref{E-Ldiscreto}), if
$\displaystyle\sum_{k=0}^{N-1}\langle\mu_{\tilde{\sigma}(k)}, X^T(k,\tilde{\sigma}(k))\rangle=0,$ then
\begin{eqnarray*}\label{ExpE-Pdiscreto}
0&=&\displaystyle\sum_{k=0}^{N-1}\langle\mathcal{E}^d_{\mu}(\sigma(k-1),\sigma(k),\sigma(k+1)),X(k,\sigma(k))\rangle
\displaystyle +\sum_{k=0}^{N-1}\left(\langle\mu_{\tilde{\sigma}(k-1)},X^{V_{\beta}}(k,\tilde{\sigma}(k-1))\rangle\right.\\
&&-
\left.\langle\mu_{\tilde{\sigma}(k)},X^{V_{\beta}}(k+1,\tilde{\sigma}(k))\rangle\right).\nonumber
\end{eqnarray*}

Now, assuming that
\begin{equation}\label{condicbordediscreta}
X(k,\sigma(k))=0,\quad \forall k\neq\{1,...,N-1\},
\end{equation}
equation (\ref{ExpE-Pdiscreto}) implies that $$\displaystyle\sum_{k=1}^{N-1}\langle\mathcal{E}^d_{\mu}(\sigma(k-1),\sigma(k),\sigma(k+1)),X(k,\sigma(k))\rangle=0,$$ for all $\mathbb{Z}$-dependent section $X:\mathbb{Z}\times Q\rightarrow TQ$ satisfying (\ref{condicbordediscreta}). Therefore, a solution $\sigma:\mathbb{Z}\rightarrow Q$ of the discrete generalized variational problem must satisfy the following system of equations
\begin{equation*}
\left\{
\begin{array}{rcl}
\displaystyle\mathcal{E}^d_{\mu}(\sigma(k-1),\sigma(k),\sigma(k+1))&=& 0 \\
\pi_{T(Q\times Q)}(\mu(k))&=& \tilde{\sigma}(k),
\end{array}
\right.
\end{equation*}
for all $1\leq k\leq N-1.$

In particular, $Im(\tilde{\sigma}(k))\subseteq C_d=\displaystyle\pi_{T(Q\times Q)}(\Sigma_d).$

If we assume that $\Sigma_d$ is determined by the vanishing of constraints $\Phi^{\alpha}=0$ in $T^*(Q\times Q),$ a curve $\sigma:\mathbb{Z}\rightarrow Q$ is a solution of the discrete generalized variational problem providing there exists a 1-form $\mu=(\mu_1,\mu_2)\in\Lambda^1(Q\times Q)$ along $\tilde{\sigma}$ such that

\begin{equation*}
\left\{
\begin{array}{rcl}
\displaystyle(\mu_1)_{(\sigma(k),\sigma(k+1))}+(\mu_2)_{(\sigma(k-1),\sigma(k))} &=& 0,\ \text{for}\ 1\leq k\leq N-1, \\
\\
\Phi^{\alpha}(\sigma(k),\sigma(k+1),(\mu_1)_{(\sigma(k),\sigma(k+1))},(\mu_2)_{(\sigma(k),\sigma(k+1))})&=& 0,\ \text{for}\ 0\leq k\leq N-1.
\end{array}
\right.
\end{equation*}

\subsection{Lagrangian mechanics}
If we have a discrete Lagrangian $L_d:Q\times Q\rightarrow\mathbb{R},$ we can consider $\Sigma_d=Im(\mathrm{d}L_d)\subseteq T^*(Q\times Q)$ and apply the discrete  generalized variational calculus. Hence, we obtain that a solution $\sigma:\mathbb{Z}\rightarrow Q$ satisfies the well-known discrete Euler-Lagrange equations (see \cite{marsden-west})
$$D_1L_d(\sigma(k),\sigma(k+1))+D_2L_d(\sigma(k-1),\sigma(k))=0,$$
for $1\leq k\leq N-1.$

\subsection{Discrete constrained variational calculus}

A \emph{discrete constrained variational problem} is defined by a pair $(C_d,l_d)$ where $C_d$ is a submanifold of $Q\times Q$ with inclusion $i_{C_d}: C_d\rightarrow Q\times Q,$ and $i_d:C_d\hookrightarrow\mathbb{R}$ is a smooth function. Now, we construct  the submanifold
\begin{multline*}
\Sigma_{l_d}=\displaystyle\left\{\mu\in T^*(Q\times Q):\pi_{T(Q\times Q)}(\mu)\in C_d\quad\text{and}\quad\langle\mu,v\rangle=\langle\mathrm{d}l_d,v\rangle,\right.\\
\left.\text{for all}\quad v\in TC_d\subseteq T(Q\times Q)
\quad\text{such that}\quad\tau_{T(Q\times Q)}(v)=\pi_{T(Q\times Q)}(\mu)\right\}\; .
\end{multline*}
In other words,
$$\Sigma_{l_d}=\left\{\mu\in T^*(Q\times Q):i^*_{C_d}\mu=\mathrm{d}l_d\right\}=\left(\mathrm{d}L_d+\nu^*(C_d)\right)\mid_{C_d},$$
where $L_d:Q\times Q\rightarrow\mathbb{R}$ is an arbitrary extension of $l_d$ to $Q\times Q$ and $\nu^*(C_d)$ is the conormal bundle.

Therefore, a solution of the discrete generalized variational problem corresponding to $\Sigma_{l_d}$ is a pair $(\sigma,\nu)$ with $\sigma:\mathbb{Z}\rightarrow Q$ and $\displaystyle\nu:\mathbb{Z}\rightarrow \nu^*(C_d)\mid_{(\sigma(k),\sigma(k+1))},$ given by the following system of difference equations:

\begin{equation}\label{DiscConstSystem}
\left\{
\begin{array}{rcl}
\displaystyle\mathcal{E}^d_{\mathrm{d}L_d+\nu}(\sigma(k-1),\sigma(k),\sigma(k+1)) &=& 0\\
(\sigma(k),\sigma(k+1))&\in & C_d.
\end{array}
\right.
\end{equation}

If $C_d$ is determined by the vanishing of constraints $\Phi^{\alpha}(q,\tilde{q})=0,$ $1\leq\alpha\leq m,$ being $m=2\dim Q-\dim C_d,$ then $\nu^*(C_d)\mid_{(\sigma(k),\sigma(k+1))}=\hbox{span}\left\{\mathrm{d}\Phi^{\alpha}(\sigma(k),\sigma(k+1))\right\}$ and $\nu(k)=\displaystyle\lambda_{\alpha}^k\mathrm{d}\Phi^{\alpha}(\sigma(k),\sigma(k+1)).$ Hence, equations (\ref{DiscConstSystem}) can be rewritten as follows (see \cite{stern})
\begin{equation}
\left\{
\begin{array}{rcl}
\displaystyle D_1(L_d+\lambda_{\alpha}^{k+1}\Phi^{\alpha})(\sigma(k),\sigma(k+1))+
D_2(L_d+\lambda_{\alpha}^{k}\Phi^{\alpha})(\sigma(k-1),\sigma(k))&=& 0\\
\text{for}\quad 1\leq k\leq N-1&&\\
\Phi^{\alpha}(\sigma(k),\sigma(k+1))= 0,\quad\text{for}\quad 0\leq k\leq N-1.&&
\end{array}
\right.
\end{equation}

\subsection{Discrete optimal control theory.}\label{DiscOptControl}

A discrete optimal control problem is specified  by a set $(U,Q,\Gamma_d,L_d)$ where $\tau_{U,Q}:U\rightarrow Q$ is a control bundle and $\Gamma_d:U\rightarrow Q\times Q$ is such that $\alpha\circ\Gamma_d=\tau_{U,Q},$ being $\alpha:Q\times Q \rightarrow Q$  the projection onto the first factor and $L_d: U\rightarrow {\mathbb R}$ is a discrete cost function. If $u_q\in U,$ then $\Gamma_d(u_q)=(q,\tilde{\Gamma}_d(u_q)).$ Taking coordinates $(q^i,u^a)$ in $U,$ we have that $\Gamma_d(q^i,u^a)=(q^i,\Gamma_{d}^i(q,u))$, that is,
$\tilde{\Gamma}_d(u_q)=(\Gamma^i(q,u))$.

As in the continuous case (see Definition \ref{discretemechanicsSolution}), the discrete optimal control solution can be obtained from the following subset of $T^*(Q\times Q):$
$$\Sigma_d=\left\{\mu\in T^*(Q\times Q):\Gamma^*_d\mu=\mathrm{d}L_d\right\}.$$

Locally, considering coordinates $(q^i,\tilde{q}^i)$ in $T^*(Q\times Q)\simeq T^*Q\times T^*Q$ we can write $\mu=(\mu_1)_i\mathrm{d}q^i+(\mu_2)_i\mathrm{d}\tilde{q}^i,$ and $\mu\in\Sigma_d$ implies that $$\Gamma_d^*\mu=\left((\mu_1)_j+(\mu_2)_i\frac{\partial\Gamma_d^i}{\partial q^j}\right)\mathrm{d}q^j+(\mu_2)_i\frac{\partial\Gamma^i_d}{\partial u^a}\mathrm{d}u^a=\frac{\partial L_d}{\partial q^j}\mathrm{d}q^j+\frac{\partial L_d}{\partial u^a}\mathrm{d}u^a=\mathrm{d}L_d.$$

Then, a solution curve $\sigma:\mathbb{Z}\rightarrow Q$ is such that there exists a curve  $\mu:\mathbb{Z}\rightarrow\Sigma_d$ given by $\mu(k)=(\mu_1(k),\mu_2(k)),$ verifying the following system of equations:
\begin{equation}
\left\{
\begin{array}{rcl}
\displaystyle (\mu_1)_j(k)+(\mu_2)_i(k)\frac{\partial\Gamma^i_d}{\partial q^j}(\sigma(k),u(k)) &=&\displaystyle\frac{\partial L_d}{\partial q^j}(\sigma(k),u(k))\\
\\
\displaystyle (\mu_2)_i(k)\frac{\partial\Gamma^i_d}{\partial u^a}(\sigma(k),u(k))&=& \displaystyle\frac{\partial L_d}{\partial u^a}(\sigma(k),u(k))\\
\\
\sigma(k+1)&=&\tilde{\Gamma}_d(\sigma(k),u(k))\\
\\
\mu_1(k+1) &=& -\mu_2(k),
\end{array}
\right.
\end{equation}
where the three  first equations are equivalent to  $\mu\in\Sigma_d,$ and the last one is equivalent to $\mathcal{E}_{\mu}^d(k,k+1)=0.$

Using the last equation in the two first of them, we obtain
\begin{equation}
\left\{
\begin{array}{rcl}
\displaystyle (\mu_1)_j(k)-(\mu_1)_i(k+1)\frac{\partial\Gamma^i_d}{\partial q^j}(\sigma(k),u(k)) &=&\displaystyle\frac{\partial L_d}{\partial q^j}(\sigma(k),u(k))\\
\\
\displaystyle -(\mu_1)_i(k+1)\frac{\partial\Gamma^i_d}{\partial u^a}(\sigma(k),u(k))&=& \displaystyle\frac{\partial L_d}{\partial u^a}(\sigma(k),u(k))\\
\\
\sigma(k+1)&=&\tilde{\Gamma}_d(\sigma(k),u(k)).
\end{array}
\right.
\end{equation}
And, if we define locally $H(q,\mu_1,u)=(\mu_1)_i\Gamma^i_d(q,u)+L(q,u),$ we obtain the following system
\begin{equation}
\left\{
\begin{array}{rcl}
\displaystyle (\mu_1)_j(k)&=&\displaystyle \frac{\partial H}{\partial q^j}(\sigma(k),\mu_1(k+1),u(k))\\
\\
0 &=& \displaystyle \frac{\partial H}{\partial u^a}(\sigma(k),\mu_1(k+1),u(k)) \\
\\
\sigma(k+1)&=&\displaystyle\frac{\partial H}{\partial \mu_1}(\sigma(k),\mu_1(k+1),u(k)).
\end{array}
\right.
\end{equation}
These are the discrete optimal control equations in this context (see \cite{lewisOptControl}).

\section{\textbf{\textsc{Discrete generalized variational calculus on Lie groupoids}}}\label{sec-groupoids}

As said in the previous section, the cartesian product $Q\times Q$ plays the role of the tangent bundle $TQ$ in the discrete setting. The geometric
relation between both spaces is expressed saying  that $Q\times Q$ has a groupoid structure being $TQ$ its associated Lie algebroid.
The purpose of this section is to describe a version of discrete generalized variational calculus adapted to general Lie groupoids  covering interesting
cases of discrete reduced dynamics (see \cite{weinstein96,MMM06Grupoides,stern}).

\subsection{Lie groupoids}

\begin{definition}
A \textbf{Lie groupoid}, denoted by $G\rightrightarrows Q$, consists of two differentiable manifolds  $G$ and $Q,$ and the
  following differentiable maps (the structural maps).
\begin{enumerate}
\item A pair of submersions: the  \textbf{source map} $ \alpha \colon G \rightarrow Q $ and the \textbf{target map} \newline $ \beta \colon G \rightarrow Q $.
\item An associative \textbf{multiplication map} $ m \colon G _2 \rightarrow G$, where
\[
  G _2 = \left\{ \left( g , h
    \right) \in G \times G \;\middle\vert\; \beta (g) = \alpha (h)
  \right\}
\]
is called the set of \textbf{composable pairs}, such that $ m(g,h) = gh. $
\item An \textbf{identity section} $ \epsilon \colon Q \rightarrow G $ of $\alpha$ and $\beta$, such
  that for all $ g \in G $,
\[
  \epsilon \left( \alpha (g) \right) g = g = g\,\epsilon\left( \beta
    (g) \right).
\]
\item An \textbf{inversion map} $ i \colon G \rightarrow G $, $ g \mapsto
  g ^{-1} $, such that for all $ g \in G $,
\[
  g g ^{-1} = \epsilon \left( \alpha (g) \right) , \qquad g ^{-1} g =
  \epsilon \left( \beta (g) \right) .
\]
\end{enumerate}
\end{definition}

Next, we will introduce the notion of a left (right) translation
by an element of a Lie groupoid. Given a groupoid $ G \rightrightarrows Q $ and an element $ g \in G
  $, we define the \textbf{left translation} $ \ell _g \colon \alpha ^{-1}
  \left( \beta (g) \right) \rightarrow \alpha ^{-1} \left( \alpha (g)
  \right) $ and the \textbf{right translation} $ r _g \colon \beta ^{-1}
  \left( \alpha (g) \right) \rightarrow \beta ^{-1} \left( \beta (g)
  \right) $ by $g$ to be
\[
  \ell _g (h) = g h , \qquad r _g (h) = h g .
\]

Analogously to the case of Lie groups, one may introduce the notion of left (right)-invariant vector field in a Lie groupoid from these translations.
Given a Lie groupoid $ G \rightrightarrows Q $, a vector field $ \xi
  \in \mathfrak{X} (G) $ is \textbf{left-invariant} if $\xi$ is $\alpha$-vertical (i.e., $T\alpha(\xi)=0$)
  and $ \left( T _h \ell_g \right) \left( \xi (h) \right) = \xi \left( g h \right), $ for
  all $(g, h) \in G _2 $.  Similarly, $\xi$ is \textbf{right-invariant}
  if $\xi$ is $\beta$-vertical  (that is, $T\beta(\xi)=0)$ and
  $ \left( T _h r _g \right) \left(  \xi (h) \right) = \xi \left( h g
  \right), $ for all $ (h, g) \in G _2 $.

It is well known that there always exists a Lie algebroid associated to a Lie groupoid (again analogously to the Lie group case).
We consider the vector bundle $\tau_{_{AG}}: AG \to Q$, whose fiber
at a point $x \in Q$ is $(AG)_{x}=\mbox{ker}\,(T_{\epsilon(x)}\alpha)$. It is easy to prove that there
exists a bijection between the space $\Gamma(\tau_{AG})$ and the set of
left (right)-invariant vector fields on $G$.
If $X$ is a section of $\tau_{_{AG}}: AG \to Q$, the corresponding
left (right)-invariant vector field on $G$
will be denoted $\overleftarrow{X}$ (respectively, $\overrightarrow{X}$), where
\begin{align}\label{defCamposInvar}
\overleftarrow{X}(g) &= (T_{\epsilon(\beta(g))}\ell_{g})(X(\beta(g))),\\
\overrightarrow{X}(g) &= -(T_{\epsilon(\alpha(g))}r_{g})((T_{\epsilon
(\alpha(g))}i)( X(\alpha(g)))),\nonumber
\end{align}
for $g \in G$. Using the above facts, we may introduce a {\bf Lie
algebroid structure} $(\lcf\cdot , \cdot\rcf, \rho)$ on
$AG$, which is defined by
\[
\overleftarrow{\lcf X, Y\rcf} = [\overleftarrow{X}, \overleftarrow{Y}], \makebox[.3cm]{}
\rho(X)(x) = (T_{\epsilon(x)}\beta)(X(x)),
\]
for $X, Y \in \Gamma(\tau_{AG})$ and $x \in Q$. Note that
\[
\overrightarrow{\lcf X, Y\rcf} = -[\overrightarrow{X}, \overrightarrow{Y}], \quad\text{and}\quad [\overrightarrow{X}, \overleftarrow{Y}] = 0,
\]
(for more details, see \cite{mackenzie87}).

\subsection{Local expressions of structural maps}

In order to obtain local expressions of the equations derived from generalized variational calculus on Lie groupoids, we need to consider some notions (see \cite{MMMLocalDiscrete} for details).  We shall begin with the notion of symmetric neighborhood.

\begin{definition}
   An open set $\mathcal{W}$ is said to be a \textbf{symmetric neighborhood} associated to an open subset $\mathcal{U}$ of a Lie groupoid $G$ if given $q_0\in Q$ a point such that $\epsilon(q_0)\in\mathcal{U},$ there exists an open subset $\mathcal{W}\subset\mathcal{U}$ of $G$ with $\epsilon(q_0)\in\mathcal{W}$ and such that
   \begin{enumerate}
      \item $\epsilon(\alpha(\mathcal{W}))\subset\mathcal{W}$ and $\epsilon(\beta(\mathcal{W}))\subset\mathcal{W},$
      \item $i(\mathcal{W})=\mathcal{W},$ and
      \item $m((\mathcal{W}\times\mathcal{W})\cap G_2)\subset\mathcal{U}.$
   \end{enumerate}
\end{definition}

\vspace{.3cm}
We consider a point
$q_0\in Q$ and a local coordinate system $(q,v)$, defined in a
neighborhood $\mathcal{U}\subset G$ of $\epsilon(q_0)$, adapted to the
fibration $\alpha:G\rightarrow M$, that is, if the coordinates of
$g\in\mathcal{U}$ are $(q^i,v^A),$ then the coordinates of
$\alpha(g)\in Q$ are $(q^i)$. We assume that the
identities correspond to elements with coordinates $(q,0)$. The
target map $\beta$ defines a local function $\textbf{b}$ as follows:
if the coordinates of $g$ are $(q,v)$, then  the coordinates of
$\beta(g)$ are $\textbf{b}(q,v)$. Note that $\textbf{b}(q,0) = q$.
Two elements $g$ and $h$ with coordinates $(q,v)$ and $(\tilde{q},\tilde{v})$ are
composable if and only if $\tilde{q}=\textbf{b}(q,v)$. Hence, local coordinates
for $G_2$ are given by $(q,v,\tilde{v})$.

 Next we consider a
symmetric neighborhood $\mathcal{W}$ associated to $q_0$ and $\mathcal{U}$. If
two elements $g, h\in\mathcal{W}$ with coordinates $(q,v)$ and $(\tilde{q},\tilde{v}),$
respectively, are composable, then the product
$gh$ has coordinates $(q,\textbf{p}(q,v,\tilde{v}))$ for some smooth function
$\textbf{p}$. We will write \newline
$
(q,v)\cdot(\tilde{q},\tilde{v})=(q,\textbf{p}(q,v,\tilde{v})).
$

We can define the following functions in terms of $\textbf{b}(q,v)$ and $\textbf{p}(q,v,\tilde{v})$,
\begin{equation}\label{rho-L-R}
\rho^i_A(q)=\frac{\partial\textbf{b}^i}{\partial v^A}(q,0),\quad \text{L}^A_B(q,v)=\frac{\partial\textbf{p}^A}{\partial \tilde{v}^B}(q,v,0)\quad\text{and}\quad
\text{R}^A_B(q,\tilde{v})=\frac{\partial\textbf{p}^A}{\partial v^B}(q,0,\tilde{v}).
\end{equation}

We will also take into account that
\begin{equation}
\label{pd and pd2 of p}
\begin{aligned}
&\frac{\partial{\textbf{p}^A}}{\partial v^B}(q,v,0)=\delta^A_B,
\qquad\qquad&&\frac{\partial{^2\textbf{p}^A}}{\partial v^B\partial v^C}(q,v,0)=0,\\
&\frac{\partial{\textbf{p}^A}}{\partial \tilde{v}^B}(q,0,\tilde{v})=\delta^A_B,
&&\frac{\partial{^2\textbf{p}^A}}{\partial \tilde{v}^B\partial v^C}(q,0,\tilde{v})=0.\\
\end{aligned}
\end{equation}
%and that the only relevant second order derivatives are given by
%\begin{equation}
%\label{local.structure.functions}
%\mathcal{C}^A_{BC}(q) \equiv \frac{\partial{^2\textbf{p}^A}}{\partial v^B\partial
%\tilde{v}^C}(q,0,0) -\frac{\partial{^2\textbf{p}^A}}{\partial \tilde{v}^B\partial v^C}(q,0,0).
%\end{equation}
%{}From the definition of $\text{L}^A_B$ and $\text{R}^A_B$ it
%follows that
%\begin{equation}
%\begin{aligned}
%\mathcal{C}^A_{BC}(q)
%=\frac{\partial{\text{L}^A_C}}{\partial v^B}(q,0)-\frac{\partial{\text{L}^A_B}}{\partial v^C}(q,0)=\frac{\partial{\text{R}^A_B}}{\partial \tilde{v}^C}(q,0)-\frac{\partial{\text{R}^A_C}}{\partial \tilde{v}^B}(q,0).
%\end{aligned}
%\end{equation}

\paragraph{\textbf{Invariant vector fields}}
If $g_0\in \mathcal{W}\subset G$ has coordinates $(q_0,v_0)$,
then the elements on the $\alpha$-fiber $\alpha^{-1}(\beta(g_0))$
have coordinates of the form $(\textbf{b}(q_0,v_0),\tilde{v})$, and the
coordinates of $l_{g_0}(g)$ are $(q_0,\textbf{p}(q_0,v_0,\tilde{v}))$. We will
write
\begin{equation}
l_{(q_0,v_0)}(\textbf{b}(q_0,v_0),\tilde{v})=(q_0,\textbf{p}(q_0,v_0,\tilde{v})).
\end{equation}

Similarly, for $h_0=(q_0,v_0)\in \mathcal{W}\subset G,$ we will write
\begin{equation}
r_{(q_0,v_0)}(q,v)=(q,\textbf{p}(q,v,v_0)).
\end{equation}

A left-invariant vector field has the form
$\overleftarrow{X}(g)=T_{\epsilon(\beta(g))}l_g(w),$ for $w\in \ker
T_{\epsilon(\beta(g))}\alpha$. To obtain a local basis of
left-invariant vector fields, we can take the local coordinate
basis $\displaystyle e_A={\frac{\partial}{\partial v^A}}\Big|_{\epsilon(\beta(g))}$ of $\ker
T_{\epsilon(\beta(g))}\alpha$. Thus, for $g\in G$ with coordinates
$(q,v)$, we have
\begin{equation*}
\label{left-alpha}
\overleftarrow{e_A}(g)=T_{\epsilon(\beta(g))}l_g\left(\frac{\partial{}}{\partial v^A}\Big|_{\epsilon(\beta(g))}\right)
=\frac{\partial{\textbf{p}^B}}{\partial\tilde{v}^A}(q,v,0)\frac{\partial{}}{\partial v^B}\Big|_g
=\text{L}^B_A(q,v)\frac{\partial{}}{\partial v^B}\Big|_{(q,v)}.
\end{equation*}

Similarly, a right-invariant vector field can be written in the form
$\overrightarrow{X}(g) = T_{\epsilon(\alpha(g))}r_g(w),$ for $w\in \ker
T_{\epsilon(\alpha(g))}\beta$.
It can be proved that a basis of right invariant vector fields is given by
\begin{equation*}
\label{right-alpha}
\begin{aligned}
\overrightarrow{e_A}(g) =T_{\epsilon(\alpha(g))}r_g\left(
   -\rho^i_A\frac{\partial{}}{\partial q^i}\Big|_{\epsilon(\alpha(g))}
   +\frac{\partial{}}{\partial v^A}\Big|_{\epsilon(\alpha(g))}\right)
=-\rho^i_A(q)\frac{\partial{}}{\partial q^i}\Big|_{g}
 +\text{R}^B_A(q,v)\frac{\partial{}}{\partial v^B}\Big|_{g},
\end{aligned}
\end{equation*}
where as before $(q,v)$ are the coordinates for $g\in G$ (see \cite{MMMLocalDiscrete} for details).

\subsection{Discrete Euler-Lagrange operator}

\vspace{.3cm}
As in section \ref{DGVCenQxQ}, we need to introduce the notion of lifts of sections of the associated Lie algebroid.

   If $X\in\Gamma(\tau_{AG}),$ we define its \textbf{complete lift} $X^c\in\mathfrak{X}(G)$ as
    $$X^c(g)=\overleftarrow{X}(g)-\overrightarrow{X}(g)\in T_gG.$$
Also, we have two notions of \textbf{vertical lift} $X^{v_{\alpha}}$ and $X^{v_{\beta}}$ of $X\in\Gamma(\tau_{AG})$ given by
$$X^{v_{\alpha}}(g)=-\overrightarrow{X}(g)\quad\text{and}\quad X^{v_{\beta}}(g)=\overleftarrow{X}(g).$$

For time dependent sections $X:\mathbb{Z}\times Q\rightarrow AG$ such that $X(k,q)\in A_qG,$ we denote $X_k\in\Gamma(AG)$ the sections given by
$X_k(q)=X(k,q),$ and then we define the \textbf{tangent lift} $X^T:\mathbb{Z}\times G\rightarrow TG$ as follows
$$X^T(k,g)=\overleftarrow{X}_{k+1}(g)-\overrightarrow{X}_k(g).$$
In the same way as in previous sections, we define $$X^{V_{\alpha}}(k,g)=X_k^{v_{\alpha}}(g)=-\overrightarrow{X}_k(g),\quad\text{and}\quad X^{V_{\beta}}(k,g)=X_k^{v_{\beta}}(g)=\overleftarrow{X}_k(g).$$

\vspace{.3cm}
\begin{definition}
 The \textbf{discrete Euler-Lagrange operator} associated with a 1-form $\mu\in\bigwedge^1G=\Gamma(\pi_{TG})$ is the mapping $\mathcal{E}_{\mu}^d:G_2\rightarrow A^*G$ defined by $$\langle\mathcal{E}_{\mu}^d(g,h),X(\beta(g))\rangle=\langle\mu_h,X^{v_{\alpha}}(h)\rangle-
 \langle\mu_g,X^{v_{\alpha}}(g)\rangle+\langle\mu_g,X^c(g)\rangle,$$
 where $\beta(g)=\alpha(h)$ since $(g,h)\in G_2.$

 Alternatively, $$\langle\mathcal{E}_{\mu}^d(g,h),X(\beta(g))\rangle=-[\langle\mu_h,X^{v_{\beta}}(h)\rangle-
 \langle\mu_g,X^{v_{\beta}}(g)\rangle-\langle\mu_h,X^c(h)\rangle],$$ for all $X\in\Gamma(\tau_{AG}).$
\end{definition}

Another useful expression is the following.
\begin{equation}\label{usefulEmu}
\langle\mathcal{E}_{\mu}^d(g,h),X(\beta(g))\rangle=\left\langle\mu_g,\overleftarrow{X}(g)\right
\rangle-\left\langle\mu_h,\overrightarrow{X}(h)\right\rangle.
\end{equation}

Therefore, using the definitions of $X^{v_{\alpha}},$ or $X^{v_{\beta}},$ and $X^c,$ we have that \begin{equation}\label{OperadorEuler-Lagrange}
\mathcal{E}^d_{\mu}(g,h)=l_g^*\mu_g+(r_h\circ i)^*\mu_h.
\end{equation}

\vspace{.3cm}
If $X:\mathbb{Z}\times Q\rightarrow AG$ with $X(k,q)\in A_qG,$ for all $k\in\mathbb{Z}$ and for all $q\in Q,$ then the Euler-Lagrange operator is
given by
\begin{eqnarray}\label{E-LdiscretoGrupoide}
\langle\mathcal{E}_{\mu}^d(g,h),X(k,\beta(g))\rangle &=&
\langle\mu_{h},X^{V_{\alpha}}(k,h)\rangle-
\langle\mu_{g},X^{V_{\alpha}}(k-1,g)\rangle+ %\nonumber \\
\langle\mu_{g},X^T(k-1,g)\rangle= \nonumber \\
& &-(\langle\mu_{h},X^{V_{\beta}}(k+1,h)\rangle-
\langle\mu_{g},X^{V_{\beta}}(k,g)\rangle-  \langle\mu_{h},X^T(k,h)\rangle). \nonumber
\end{eqnarray}
\vspace{.5cm}

\textbf{Local expressions.} Locally, if $(g,h)\in G_2$ is a composable pair that both $g$ and $h$ are on the same symmetric neighborhood
$\mathcal{W}$, with coordinates $(q,v)=(q^i,v^A)$ for $g$ and $(\tilde{q},\tilde{v})=(\tilde{q}^i,\tilde{v}^A)$ for $h$ in the Lie groupoid $G,$ we
can write $\mu=(\mu_1)_i(q,v)\mathrm{d}q^i+(\mu_2)_A(q,v)\mathrm{d}v^A,$ and using expression (\ref{usefulEmu}) to compute the discrete Euler-Lagrange
operator in the base $\{e_A\}$ of $\Gamma(\tau_{AG}),$ we get
\begin{eqnarray}\label{eclocal}
\left\langle\mathcal{E}_{\mu}^d(q,v,\tilde{v}),e_A(\beta(q,v))\right\rangle&=&
(\mu_1)_i(\tilde{q},\tilde{v})\rho^i_A(\tilde{q})+(\mu_2)_B(q,v)\text{L}^B_A(q,v)\\
&&-(\mu_2)_B(\tilde{q},\tilde{v})\text{R}^B_A(\tilde{q},\tilde{v}),\nonumber
\end{eqnarray}
where $\tilde{q}=\textbf{b}(q,v)$.
%\vspace{.2cm}
%Then, we obtain the following system of equations

%\begin{eqnarray}\label{E-LLieGroupoid}
  % (\mu_1)_i\rho^i_A(\tilde{q})+(\mu_2)_BL^B_A(q,v)
%-(\mu_2)_BR^B_A(\tilde{q},\tilde{v})&=&0\\
%\tilde{v}&=&\textbf{b}(q,v)\nonumber
%\end{eqnarray}

\subsection{Discrete mechanics on Lie groupoids}

\begin{definition}
   A \textbf{discrete generalized variational problem} on a Lie groupoid $G$ is determined by a submanifold $\Sigma_d\subseteq T^*G.$
\end{definition}

 Given a curve $\gamma:\mathbb{Z}\rightarrow G,$ we say that $\gamma$ is \textbf{composable} if $\beta(\gamma(k))=\alpha(\gamma(k+1)).$

\begin{definition}
   A \textbf{solution} of a discrete generalized variational problem determined by  $\Sigma_d\subseteq T^*G$ is a composable curve
   $\gamma:\mathbb{Z}\rightarrow G$ such that there exists $\mu:\mathbb{Z}\rightarrow\Sigma_d$ verifying $\pi_{TG}(\mu(k))=\gamma(k),$ and for all
   $\mathbb{Z}$-dependent section $X:\mathbb{Z}\times Q\rightarrow AG,$ we have that $$\sum_{k=0}^{N-1}\langle\mu(\gamma(k)),X^T(k,\gamma(k))\rangle=0.$$
\end{definition}

\vspace{.2cm}
Since (\ref{E-LdiscretoGrupoide}) holds, if
$\displaystyle\sum_{k=0}^{N-1}\left\langle\mu(\gamma(k)),X^T(k,\gamma(k))\right\rangle=0,$ then
\begin{eqnarray}\label{ExpE-PdiscretoGrupoide}
&\displaystyle\sum_{k=0}^{N-1}&\left\langle\mathcal{E}^d_{\mu}\left(\gamma(k),\gamma(k+1)\right),
X(k,\beta(\gamma(k)))\right\rangle=\\
&\displaystyle\sum_{k=0}^{N-1}&\left(\left\langle\mu_{\gamma(k)},X^{V_{\beta}}(k,\gamma(k))\right\rangle-
\left\langle\mu_{\gamma(k+1)},X^{V_{\beta}}(k+1,\gamma(k+1))\right\rangle\right).\nonumber
\end{eqnarray}

\vspace{.3cm}
Now, if we assume the boundary conditions
\begin{equation}\label{condicbordediscretaGrupoide}
X(k,\beta(\gamma(k)))=0,\quad \forall k \neq\{1,...,N-1\},
\end{equation}
equation (\ref{ExpE-PdiscretoGrupoide}) and expression (\ref{defCamposInvar}) imply that
$$\displaystyle\sum_{k=0}^{N-1}\left\langle\mathcal{E}^d_{\mu}\left(\gamma(k),\gamma(k+1)),
X(k,\beta(\gamma(k))\right)\right\rangle=0, $$ for all $\mathbb{Z}$-dependent section $X:\mathbb{Z}\times Q\rightarrow AG$ satisfying
(\ref{condicbordediscretaGrupoide}). Therefore, a solution $\gamma:\mathbb{Z}\rightarrow G$ of the discrete generalized variational problem on a Lie
groupoid $G$ must satisfy the following system of equations
\begin{equation}\label{main-groupoids}
\left\{
\begin{array}{rcl}
\displaystyle\mathcal{E}^d_{\mu}\left(\gamma(k),\gamma(k+1)\right)&=& 0 \\
\pi_{TG}(\mu(k))&=& \gamma(k)\\
\beta(\gamma(k))&=&\alpha(\gamma(k+1))\; ,
\end{array}
\right.
\end{equation}
for all $0\leq k\leq N-1.$

\subsection{Lagrangian mechanics}

If we have $L_d:G\rightarrow\mathbb{R},$ we can take $\Sigma_d=Im(\mathrm{d}L_d)\subseteq T^*G$ and applying discrete generalized variational calculus to
obtain that the discrete Euler-Lagrange equation is
\begin{eqnarray}\label{E-LDiscdL}
\mathcal{E}^d_{\mathrm{d}L_d}(\gamma(k),\gamma(k+1))&=&\displaystyle l_{\gamma(k)}^*\mathrm{d}L_d(\gamma(k))+(r_{\gamma(k+1)}\circ i)^*\mathrm{d}L_d(\gamma(k+1))\\
&=&\mathrm{d}(L_d\circ l_{\gamma(k)}+L_d\circ r_{\gamma(k+1)}\circ i)(\epsilon(\beta(\gamma(k))))=0\nonumber
\end{eqnarray}
or, in other words, $$\langle\mathrm{d}L_d(\gamma(k)),\overleftarrow{X}(\gamma(k))\rangle-
\langle\mathrm{d}L_d(\gamma(k+1)),\overrightarrow{X}(\gamma(k+1))\rangle=0,$$
for all $X\in\Gamma(\tau_{AG})$ and $0\leq k\leq N-1.$ %(see \ref{E-LDiscdL}).

Equivalently, we can write $$\overleftarrow{X}(\gamma(k))(L_d)-\overrightarrow{X}(\gamma(k+1))(L_d)=0,$$
for $0\leq k\leq N-1.$
These equations coincide with the ones given in \cite{MMM06Grupoides}.

\vspace{.3cm}
Locally, if $\gamma(k)=(q_k, v_k),$ these equations are (see \ref{eclocal}):
\begin{eqnarray*}
0&=&\frac{\partial L_d}{\partial q^i}(q_{k+1}, v_{k+1})\rho^i_A(q_{k+1})+\frac{\partial L_d}{\partial v^B}(q_{k}, v_{k})\text{L}^B_A(q_{k}, v_{k})\\
&&-\frac{\partial L_d}{\partial v^B}(q_{k+1}, v_{k+1})\text{R}^B_A(q_{k+1}, v_{k+1})\\
q_{k+1}&=&\textbf{b}(q_k,v_k).
\end{eqnarray*}

\subsubsection{Discrete Euler-Poincar\'e equations} If we consider a Lie group $G$ as a groupoid over the identity element $e\in G,$ the structural maps
are $\alpha(g)=e,$ $\beta(g)=e,$ $\epsilon(e)=e,$ $i(g)=g^{-1},$ $m(g,h)=gh,$ for all $g,h\in G.$ The Lie algebroid associated with $G$ is the Lie algebra
$\mathfrak{g}=T_eG$ of $G$ and, given $\xi\in\mathfrak{g},$ we have that the left and right-invariant vector fields are:
$\overleftarrow{\xi}(g)=(T_el_g)(\xi)$ and $\overrightarrow{\xi}(g)=(T_er_g)(\xi),$ for $g\in G.$ Hence, given a Lagrangian $L:G\rightarrow\mathbb{R},$
the \emph{discrete Euler-Poincar\'e equations} are
$$\mathcal{E}_{\mathrm{d}L_d}^d=\mathrm{d}(L_d\circ l_{\gamma(k)}-L_d\circ r_{\gamma(k+1)})(e)=0$$ (see equation \ref{E-LDiscdL}), or
$(l^*_{\gamma(k)}\mathrm{d}L)(e)=(r^*_{\gamma(k+1)}\mathrm{d}L)(e).$
Equivalently, we have $$(T_el_{\gamma(k)})(\xi)(L)-(T_er_{\gamma(k+1)})(\xi)(L)=0,$$ for all $\xi\in\mathfrak{g}.$
These equations coincide with the ones obtained in \cite{MMM06Grupoides}. If we denote $\mu_k=(r^*_{\gamma(k)}{d}L)(e),$ the discrete Euler-Poincar\'e
equations are written as $\mu_{k+1}=\text{Ad}^*_{\gamma(k)}\mu_k,$ where $\text{Ad}:G\times\mathfrak{g}\rightarrow\mathfrak{g}$ is the adjoint
action of $G$ on $\mathfrak{g}.$ These equations are known as the \emph{discrete Lie-Poisson equations}.

\subsection{Discrete constrained mechanics on Lie groupoids} See \cite{stern} for details about the topic.

A \emph{discrete constrained variational problem} is defined by a pair $(C_d,l_d)$ where $C_d$ is a submanifold of a Lie groupoid $G$ with inclusion
$i_{C_d}:C_d\hookrightarrow G,$ and $l_d:C_d\rightarrow\mathbb{R}$ is a function. Now, we consider the submanifold
\begin{multline*}
\Sigma_{l_d}=\displaystyle\left\{\mu\in T^*G:\pi_{TG}(\mu)\in C_d\quad\text{and}\quad\langle\mu,v\rangle=\langle\mathrm{d}l_d,v\rangle,\right.\\
\left.\text{for all}\quad v\in TC_d\subseteq TG
\quad\text{such that}\quad\tau_{TG}(v)=\pi_{TG}(\mu)\right\}.
\end{multline*}
In other words,
$$\Sigma_{l_d}=\left\{\mu\in T^*G:i^*_{C_d}\mu=\mathrm{d}l_d\right\}=\left(\mathrm{d}L_d+\nu^*(C_d)\right)\mid_{C_d},$$
where $L_d:G\rightarrow\mathbb{R}$ is an arbitrary extension of $l_d$ to $G,$ and $\nu^*(C_d)$ is the associated conormal bundle.
%defined by $\displaystyle\nu^*(C_d)=\displaystyle\left\{\nu\in T^*G\mid_{C_d}:\nu\mid_{T_{\pi_{TG}(\nu)}C_d}=0\right\}.$

Therefore, a solution of the discrete generalized variational calculus corresponding to $\Sigma_{l_d}$ is a pair $(\gamma,\nu)$ with
$\gamma:\mathbb{Z}\rightarrow Q$ and $\displaystyle\nu:\mathbb{Z}\rightarrow \nu^*(C_d)\mid_{(\gamma(k))},$ given by the following system of equations:

\begin{equation}\label{eqDiscConstSystem}
\left\{
\begin{array}{rcl}
\displaystyle\mathcal{E}^d_{\mathrm{d}L_d+\nu}(\gamma(k-1),\gamma(k)) &=& 0\\
\gamma(k)&\in & C_d\\
\beta(\gamma(k))&=&\alpha(\gamma(k+1)).
\end{array}
\right.
\end{equation}

If $C_d$ is determined by the vanishing of constraints $\Phi^{\alpha}(g)=0,$ $1\leq\alpha\leq m,$ being $m=2\dim Q-\dim C_d,$ then
$\nu^*(C_d)\mid_{(\gamma(k))}=\hbox{span}\left\{\mathrm{d}\Phi^{\alpha}(\gamma(k))\right\}$ and
$\nu(k)=\displaystyle\lambda_{\alpha}^k\mathrm{d}\Phi^{\alpha}(\gamma(k)).$ Hence, equations (\ref{eqDiscConstSystem}) can be rewritten as follows
\begin{equation*}
\left\{
\begin{array}{rcl}
\displaystyle D_1(L_d+\lambda_{\alpha}^{k+1}\Phi^{\alpha})(\gamma(k))+
D_2(L_d+\lambda_{\alpha}^{k}\Phi^{\alpha})(\gamma(k-1)) &=& 0,\ \text{for}
\ 1\leq k\leq N-1;\\
\Phi^{\alpha}(\gamma(k))&=& 0,\ \text{for}\ 0\leq k\leq N-1.
\end{array}
\right.
\end{equation*}

Locally, if $\gamma(k)=(q_{k}, v_{k}),$ these equations are:
\begin{eqnarray*}
0&=&\frac{\partial (L_d+\lambda_{\alpha}^{k+1}\Phi^{\alpha})}{\partial q^i}(q_{k}, v_{k})\rho^i_A(q_{k})+\frac{\partial( L_d+\lambda_{\alpha}^{k}\Phi^{\alpha})}{\partial v^B}(q_{k-1}, v_{k-1})\text{L}^B_A(q_{k-1}, v_{k-1})\\
&&-\frac{\partial (L_d+\lambda_{\alpha}^{k+1}\Phi^{\alpha})}{\partial v^B}(q_{k}, v_{k})\text{R}^B_A(q_{k}, v_{k})\\
0&=& \Phi^{\alpha}(q_{k}, v_{k})\\
q_{k}&=&\textbf{b}(q_{k-1},v_{k-1})\; .
\end{eqnarray*}

\subsection{Discrete optimal control theory on Lie groupoids}

A \emph{discrete optimal control problem} on a Lie groupoid $G$ is given by a set $(U,Q,\Gamma_d,L_d)$ where $\tau_{U,Q}:U\rightarrow Q$ is a control
bundle, $\Gamma_d:U\rightarrow G$ is such that $\alpha\circ\Gamma_d=\tau_{U,Q},$ being $\alpha:G \rightarrow Q$ the projection, and
$L_d: U\rightarrow {\mathbb R}$ is a discrete cost function (see \cite{MR2050148}).%That is, if $u_q\in U,$ then $\Gamma_d(u_q)=(q,\tilde{\Gamma}_d(u_q)).$

As we saw in section \ref{DiscOptControl}, the discrete optimal control solution can be obtained from the following subset of $T^*G:$
$$\Sigma_d=\left\{\mu\in T^*G\; \mid \; \Gamma^*_d\mu=\mathrm{d}L_d\right\}.$$
In local coordinates on the Lie groupoid, we obtain the following system of equations:

\begin{eqnarray*}
\frac{\partial L_d}{\partial q^i}(q_k, u_k)&=&(\mu_1)_i(q_k,u_k)+(\mu_2)_A(q_k,u_k)\frac{\partial \Gamma^A_d}{\partial q^i}(q_{k}, u_{k})\\
\frac{\partial L_d}{\partial u^A}(q_k, u_k)&=&(\mu_2)_B(q_k,u_k)\frac{\partial \Gamma^B_d}{\partial u^A}(q_{k}, u_{k})\\
%_k^A&=&\Gamma^A_d(q_{k}, u_{k})\\
q_{k+1}&=&\textbf{b}(\Gamma_d(q_k,u_k)).
\end{eqnarray*}

\section{Conclusions}

In this paper, we have introduced many of the most important equations of motion of  mechanical systems using a generalization of variational calculus
where the main in\-gre\-dient is played by a subset of the cotangent space of the velocity phase space. Cases like standard Lagrangian mechanics, nonholonomic mechanics, constrained variational calculus, hamiltonian mechanics, systems admitting a Lie group of symmetries, among others, are naturally included in this framework.
Moreover, it is possible to extend this technique to the case of discrete mechanics using a parallel construction.

In the future, we will study how the constraint algorithms work in the setting of ge\-ne\-ra\-li\-zed variational calculus, and the extension of our method to
the case of discrete nonholonomic mechanics (see \cite{iglesiaspadron}) and discrete hamiltonian systems. In our future work, we will also develop other topics such as ge\-ne\-ra\-li\-zed variational calculus both in the case of Dirac structures modeling mechanics and the theory of interconnection.

\bibliographystyle{plain}
%\bibliography{Referencias3-11-14}

\end{document}